\documentclass[MSNbibl,number,citesort,seceqn,dvips]{arxbj}
\usepackage{upgreek}
\usepackage{graphicx}

\aid{0}
\volume{19}
\issue{1}
\pubyear{2013}
\firstpage{308}
\lastpage{343}
\doi{10.3150/11-BEJ401}

\newcommand{\eqref}[1]{(\ref{#1})}
\newtheorem{Theorem}{Theorem}[section]
\newtheorem{Lemma}[Theorem]{Lemma}
\newtheorem{Proposition}[Theorem]{Proposition}
\newproclaim{Assumption}[Theorem]{Assumption}
\newremark{Remark}[Theorem]{Remark}
\def\la{\leftarrow}
\def\eqd{\stackrel{d}{=}}
\def\cinP{\stackrel{P}{\to}}
\newcommand{\sgn}{\operatorname{sgn}}
\def\E{\mathbb{E}}
\def\R{\mathbb{R}}
\def\C{\mathbb{C}}
\def\D{\mathbb{D}}
\def\F{\mathcal{F}}
\def\weak{\Rightarrow}
\begin{document}
\begin{frontmatter}

\title{Weak limits for exploratory plots in the analysis of extremes}
\runtitle{Exploratory plots in extremes}

\begin{aug}
\author[1]{\fnms{Bikramjit} \snm{Das}\corref{}\thanksref{1}\ead[label=e1]{bikram@math.ethz.ch}}%
\and
\author[2]{\fnms{Souvik} \snm{Ghosh}\thanksref{2}\ead[label=e2]{souvik@yahoo-inc.com}}
\runauthor{B. Das and S. Ghosh} 
\address[1]{RiskLab, Department of Mathematics, ETH Zurich, R\"amistrasse 101, 8092 Zurich, Switzerland.\\
\printead{e1}}
\address[2]{Yahoo! Research, 4401 Great America Parkway, Santa Clara,
CA 95054, USA.\\
\printead{e2}}
\end{aug}

\received{\smonth{9} \syear{2010}}
\revised{\smonth{8} \syear{2011}}

%
\begin{abstract}
Exploratory data analysis is often used to test the goodness-of-fit of
sample observations to specific target distributions. A few such
graphical tools have been
extensively used to detect subexponential or heavy-tailed behavior in
observed data.
In this paper we discuss asymptotic limit behavior of two such
plotting tools: the quantile--quantile plot and the mean excess plot.
The weak consistency
of these plots to fixed limit sets in an appropriate topology of $\R^2$
has been shown in Das and Resnick (\textit{Stoch. Models} \textbf{24}
(2008) 103--132) and Ghosh and Resnick
(\textit{Stochastic Process. Appl.} \textbf{120} (2010) 1492--1517). In
this paper we
find asymptotic
distributional limits for these plots when the underlying distributions
have regularly varying right-tails. As an application we construct
confidence bounds around the plots
which enable us to statistically test whether the underlying
distribution is heavy-tailed or not.
\end{abstract}

%

\end{frontmatter}

\section{Introduction}\label{secintro}
Statistical analysis of extremes in available data has been very
important in varied areas like finance (McNeil, Frey and Embrechts \cite
{mcneilfreyembrechts2005}),
telecommunication (Maulik, Resnick and Rootz{\'e}n~\cite
{maulikresnickrootzen2002}, D'Auria and Resnick \cite
{dauriaresnick2006}), hydrology (Katz, Parlange and Naveau \cite
{katzparlangenaveau2002}), environmental statistics
(Davison and Smith~\cite{davisonsmith1990}, Smith~\cite{smith2003}) and
many more. Before
analyzing features of the data using extreme value analysis, it is
imperative that
we check whether extreme-value modeling is well suited in the given
context; see Drees~\cite{drees2011}
for a recent survey of exploratory techniques for extremes in an
actuarial context. Popular exploratory techniques in this direction
have been the mean excess (ME) plots
(Davison and Smith~\cite{davisonsmith1990}) and the quantile--quantile
(QQ) plots
which are specifically tuned for heavy-tailed data (Kratz and Resnick
\cite
{kratzresnick1996}).
A distribution~$F$
is \textit{heavy-tailed} if the tail probability $(1-F)$ is regularly
varying (Resnick~\cite{resnickbook2008}, Chapter~1).
It has been shown earlier that under an assumption of heavy tails and
with proper normalizations, both plots converge in probability to fixed
closed sets (for ME plots, an additonal assumption of finiteness of the
mean of $F$ is required);
see Das and Resnick~\cite{dasresnick2008} and Ghosh and Resnick \cite
{ghoshresnick2009}. These
results corroborate the use of the QQ plot and the ME plot to test the
null hypothesis that the underlying distribution is heavy-tailed. The
proximity of the observed plot to the fixed limit set would support the
null hypothesis.


Incidentally, one data set leads to just one single plot of each kind.
A single plot is often not enough to statistically detect proximity
between the plot
and the intended fixed limit set; see the examples in Section~\ref{secdata}.
Creating appropriate confidence bounds around these plots, though, can
help us to test the null hypothesis with some degree of confidence.
In this paper we study weak limits of both kinds of plots for
heavy-tailed data and use these limits to obtain
confidence bounds around them with asymptotic coverage probabilities.
The methods used here are general and can be used to find weak limits
and confidence bounds for other plots used in the analysis of extremes.


\subsection{Plan for this paper}
We introduce the two plotting methodologies in Section~\ref{secintro}.
In Section~\ref{secprelim} we set up necessary tools to talk about
convergence of random closed sets in $\R^2$, since the QQ and ME plots
are random closed sets in $\R^2$.
In Sections~\ref{secQQplot} and~\ref{secMEplot}, we prove weak
convergence of the QQ and the ME plot under the null hypothesis that
the underlying distribution $F$ is heavy-tailed. We proceed by
expressing both plots as appropriate functionals of the tail empirical
measure and then use convergence properties of the tail empirical
measure to prove weak convergence of both plots. As an application to
obtaining these weak limits,
we construct confidence bounds with asymptotic coverage probability
for both kinds of plots in Section~\ref{secconfidence}.
Finally, in Section~\ref{secdata} we apply the results obtained in the
previous sections to simulated and real data sets
to exemplify how they perform in practice. We conclude in Section \ref
{secconclusion} along with a discussion on future directions.
\subsection{QQ plots for heavy tails}
Suppose we want to test the null hypothesis that observations from a
sample are independent and identically distributed (i.i.d.)
with some known distribution $F$. The QQ plot, which is a plot of the
empirical quantiles from the data against the distributional quantiles
of $F$, is an intuitive and popular graphical tool for detecting the
goodness-of-fit for a sample to the distribution $F$. If the true
distribution of the sample is $ F$, then the QQ plot should converge,
in an appropriate sense, to a straight line. Results involving
empirical process and quantile process convergences are available in
Shorack and Wellner~\cite{shorackwellner1986}, which can be
appropriately used to create
confidence intervals for QQ plots. The QQ plot we consider is a little
different and is specifically designed to check for distributions $F$
where $\bar F:=1-F$ is regularly varying with some index $-1/\xi$, $\xi
>0$, also denoted $\bar F \in RV_{-1/\xi}$ (Resnick \cite
{resnickbook2008}, Chapter~1).
For a sample $X_1,X_2,\ldots,X_n$, its
decreasing order statistics are denoted by $X_{(1)}\ge X_{(2)} \ge
\cdots\ge X_{ (n)}$ and the QQ plot in this context is defined by
\[
\mathcal{Q}_n = \biggl\{\biggl( -\log\frac jk, \log\frac
{X_{(j)}}{X_{(k)}} \biggr)\dvt1\le j \le k \biggr\}, \qquad k <n.
\]
Clearly we concentrate on the top $k$ quantiles of the data justified
by the fact that $\bar F \in RV_{-1/\xi}$ only provides us with
information about the right tail of the data. Under the null hypothesis
of $\bar F \in RV_{-1/\xi}$ for some $\xi>0$, Das and Resnick \cite
{dasresnick2008}
have shown convergence in probability for QQ plots in an appropriate
topology of random closed sets when the data is assumed to be an i.i.d. sample.

\subsection{ME plots}
The ME function of a random variable $X$ is defined as
%
%
\begin{equation}\label{eqmeanexcess}
M(u):= E[ X-u|X>u ],
\end{equation}
provided $EX_+<\infty,$ and is also known as the \textit{mean residual
life function}. A natural estimate of $M(u)$ is the empirical ME function
$ \hat M(u)$ defined as
%
%
\begin{equation}\label{eqempiricalmeanexcess}
\hat M(u)= \frac{ \sum_{ i=1}^{ n}(X_{ i}-u)I_{ [ X_{ i}>u]}}{ \sum_{
i=1}^{ n}I_{ [ X_{ i}>u]}},\qquad u\ge0.
\end{equation}
The ME plot is the plot of the points $ \{ (X_{ (k)},\hat M(X_{
(k)}))\dvt 1< k\le n \}$.

The ME plot is often used as a simple graphical test to check if data
conform to a generalized Pareto distribution (GPD). The GPD is an
important class of distributions and is fundamental for the
peaks-over-threshold method used in extreme value analysis (Davison and
Smith \cite
{davisonsmith1990}). The GPD is characterized by its cumulative
distribution function $ G_{ \xi,\beta}$
%
%
\begin{equation}\label{eqGPD}
G_{ \xi,\beta}(x)= \cases{
1-(1+ \xi x/ \beta)^{ -1/\xi} & \quad $\mbox{if }\xi\neq
0,$\vspace*{2pt}\cr
1-\exp(-x/\beta) & \quad $\mbox{if } \xi=0,$}
\end{equation}
where $ \beta>0$, and $ x\ge0$, when $ \xi\ge0$ and $ 0\le x\le
-\beta/\xi$, if $ \xi<0$. The parameters $ \xi$ and $ \beta$ are
referred to as the \textit{shape} and the \textit{scale} parameter,
respectively.
The GPD in the case $\xi>0$ corresponds to the classical Pareto law
with tail exponent $1/\xi$.
For a random variable $ X\sim G_{ \xi,\beta}$, we have $E(X)<\infty
$, if and only if $ \xi<1$, and in this case, the ME function of $X$
is linear in $u$.
%
%
\begin{equation}\label{eqmegdp}
M(u)= \frac{ \beta}{1-\xi}+ \frac{ \xi}{1-\xi}u,
\end{equation}
where $ 0\le u< \infty$ if $ 0\le\xi<1$ and $ 0\le u\le-\beta/\xi$
if $ \xi<0$. In fact, the linearity of the ME function characterizes
the GPD class; cf. McNeil, Frey and Embrechts \cite
{mcneilfreyembrechts2005} and Embrechts, Kl\"{u}ppelberg and Mikosch
\cite
{embrechtskluppelbergmikosch1997}. Davison and Smith \cite
{davisonsmith1990} used
this property and suggested that if the ME plot is close to a straight
line for high values of the threshold, then there is no evidence
against the use of a GPD model. See also Embrechts, Kl\"{u}ppelberg and
Mikosch \cite
{embrechtskluppelbergmikosch1997} and Hogg and Klugman~\cite
{hogg1984ld} for the
implementation of this plot in practice. Ghosh and Resnick \cite
{ghoshresnick2009}
discuss convergence in probability for the high thresholds of suitably
normalized ME plots in an appropriate topology of random closed sets
when the data is an i.i.d. sample.

The advantage of the ME plot over the QQ plot is that it works when $
-\infty<\xi< 1$, whereas the QQ plot works for $\xi>0$ only. Hence the
ME plot can be used whenever the sample is in the maximal domain of
attraction of any generalized extreme value distribution with finite
mean (Gumbel, Weibull or Fr\'echet distribution). The QQ plot is
restricted to the domain of attraction of Fr\'echet distribution only.
In this paper, though, we restrict to the case when $ \xi>0$, which is
the case of maximal domain of attraction of the Fr\'echet distribution.
The disadvantage of the ME plot is that it requires $ \xi<1$ to make
proper sense of the result,
that is, the underlying distribution should have a finite mean. Still,
limits can and have been obtained for the ME plots, even when the
distributional mean is not finite; see Ghosh and Resnick \cite
{ghoshresnick2009}.

\section{Preliminaries}\label{secprelim}
\subsection{Topology on closed sets of $\mathbb{R}^{2}$}\label{subsectop}

Since we are dealing with plots which are closed sets in $ \mathbb
{R}^{ 2}$, it is imperative to understand the topology on closed sets.
We denote the collection of all closed (compact) sets in $ \mathbb{R}^{
2}$ by $ \mathcal{F}$ ($ \mathcal{K}$, resp.). We consider a hit
and miss topology on $ \mathcal{F}$ called the Fell topology. The Fell
topology is the smallest topology containing the families $ \{ \mathcal
{F}^{ K}, K \mbox{ compact} \}$ and $ \{ \mathcal{F}_{ G}, G\mbox{
open} \}$ where, for any set $ B$,
\[
\mathcal{F}^{ B}= \{ F\in\mathcal{F}\dvt F\cap B =\varnothing\} \quad\mbox{and}\quad \mathcal{F}_{ B}=\{ F\in\mathcal{F}\dvt F\cap B \neq\varnothing
\}.
\]
Hence $ \mathcal{F}^{ B}$ and $ \mathcal{F}_{ B}$ are collections of
closed sets which miss and hit the set $ B$, respectively. This is the
reason for which such topologies are called hit and miss topologies. In
the Fell topology, a sequence of closed sets $ \{ F_{ n} \}$ converges
to $ F\in\mathcal{F}$ if and only if the following two conditions hold:
\begin{itemize}
\item$ F \in\mathcal{F}_G$ implies there exists $ N\ge1$ such that
for all $ n\ge N$, $ F_{ n} \in\mathcal{F}_G$, for any open set $G$.
%
\item$ F \in\mathcal{F}^K$ implies there exists $ N\ge1$ such that
for all $ n\ge N$, $ F_{ n}\in\mathcal{F}^K$, for any compact set~$K$.
\end{itemize}
The Fell topology on the closed sets of $ \mathbb{R}^{ 2}$ is
metrizable (Flachsmeyer~\cite{flaschmeyer1963}, Beer~\cite{beer1993a})
and we
indicate convergence in this topology of a sequence $(F_n)$ of
closed sets to a limit closed set $F$ by $F_n \to F$.
Often though, it is easier to deal with the following characterization
of convergence.
\begin{Lemma}\label{lemsetconv}
A sequence $ F_{ n}\in\mathcal{F}$ converges to $ F\in\mathcal{F}$ in
the Fell topology if and only if the following two conditions hold:
\begin{enumerate}
\item For any $ t\in F$ there exists $ t_{ n}\in F_{ n}$ such that $
t_{ n}\to t.$
\item If, for some subsequence $ (m_{ n})$, $ t_{ m_{ n}}\in F_{ m_{
n}}$ converges, then $ \lim_{ n\to\infty} t_{ m_{ n}}\in F$.
\end{enumerate}
\end{Lemma}

See Theorem 1-2-2 in Matheron~\cite{matheron1975}, page~6, for a proof
of this
lemma.

Let $\sigma_{\F}$ denote the Borel $\sigma$-algebra generated by
the Fell topology of open sets {(not to be confused with open sets in
$\R^d$)}.
A \textit{random closed set} $X\dvtx \Omega\mapsto\F$ is a measurable
mapping from $(\Omega,\mathcal{A},P^{\prime})$ to $ (\F, \sigma_{\F}
)$. Denote by $P$ the induced probability on $\sigma_{\F}$,
that is, $P=P^{\prime}\circ X^{-1}.$

Since the Fell topology is metrizable, the definition of convergence in
probability is obvious. The following result is a well-known and
helpful characterization for convergence in probability of random
variables, and it holds for random sets as well; see Theorem~6.21 in
Molchanov~\cite{molchanov2005}, page~92.

\begin{Lemma}\label{lemconvprob}
A sequence of random sets $ (F_{ n})$ in $ \mathcal{F}$ converges in
probability to a random set $ F$ if and only if for every subsequence $
(n^{ \prime})$ of $ \mathbb{Z}_{ +}$ there exists a further
subsequence $ (n^{ \prime\prime})$ of $ (n^{ \prime})$ such that $
F_{ n^{ \prime\prime}} \to F$-a.s.
\end{Lemma}

A sequence of random closed sets $(X_n)_{n\ge1} $ weakly
converges to a random closed set $X$ with distribution $P$ if the
corresponding induced
probability measures $(P_n)_{n\ge1} $ converge weakly to~$P$,
that is,
\[
P_n(\mathcal{B})=P_n^{\prime}\circ X_n^{-1}( \mathcal{B})
\quad\to \quad P(\mathcal{B})= P^{\prime}\circ X^{-1}( \mathcal{B}),
\qquad\mbox{as }  n \to\infty
\]
for each $ \mathcal{B} \in\sigma_{\mathcal{F}} $ such that
$P(\partial\mathcal{B}) =
0$. This is not always straightforward to verify from the definition.
The following characterization of weak convergence in
terms of sup-measures is very useful; cf. Vervaat~\cite{vervaat1997}. Suppose
$h\dvtx \R^d \to
\R_{+} =[0,\infty)$. For $X\subset\mathbb{R}^d$, define $h(X) = \{
h(x)\dvtx x \in X \}$, and $h^\vee$ is the sup-measure generated by $h$
defined by
\[
h^\vee(X)=\sup\{ h(x)\dvtx x \in X \}
\]
(Molchanov~\cite{molchanov2005}, Vervaat~\cite{vervaat1997}).
These definitions permit the following characterization
(Molchanov~\cite{molchanov2005}, page~87).
\begin{Lemma} \label{thmchoq}
A sequence $(X_n)_{n\ge1}$ of random closed sets converges weakly to
a random closed set $X$ if and only if $\E h^{\vee} (X_n)$
converges to $\E h^{\vee} (X)$ for every non-negative continuous
function $h\dvtx\R^d \to\R$ with a bounded support.
\end{Lemma}

We often use the following notation: for a $x \in\R$ and a
set $ A\subset\mathbb{R}^{ n}$, $xA=\{ xy\dvt y\in A \} $ and $ x+A=\{
x+y\dvt y\in A \}$. See Matheron~\cite{matheron1975} and Molchanov \cite
{molchanov2005} for
further details on the theory of random sets.

\subsection{Miscellany}\label{subsecmisc}
Throughout this paper we will take $ k:=k_{ n}$ to be a sequence
increasing to infinity such that $ k_{ n}/n\to0$. For a distribution
function $F(x)$ we write $\bar F(x):=1-F(x) $ for the tail, and the
quantile function is
\[
b(u):=F^\leftarrow\biggl(1-\frac1u\biggr)=\inf\biggl\{s\dvt F(s)\geq
1-\frac1u\biggr\}=\biggl(\frac{1}{1-F}\biggr)^\leftarrow(u).
\]
A function
$U\dvtx(0,\infty)\to\mathbb{R}_+$ is regularly varying with index
$\rho\in\mathbb{R}$, written $U\in RV_\rho$, if
\[
\lim_{t\to\infty} \frac{U(tx)}{U(t)}=x^\rho,\qquad x>0.
\]
Regular variation is discussed in several books
such as Resnick~\cite{resnickbook2007,resnickbook2008},
Seneta~\cite{seneta1976},
Geluk and de~Haan~\cite{gelukdehaan1987}, de~Haan~\cite{dehaan1970},
de~Haan and Ferreira~\cite{dehaanferreira2006},
Bingham, Goldie and Teugels~\cite{binghamgoldieteugels1987}.

We use $ M_{ +}(0,\infty]$ to denote the space of non-negative Radon
measures $ \mu$ on $
(0,\infty]$ metrized by the vague metric. Point measures are written as
a function of their points $\{x_i,
i=1,\dots,n\}$ by $\sum_{i=1}^n \delta_{x_i};$ see, for example,
Resnick~\cite{resnickbook2008}, Chapter~3.

We will use the following notations to denote different classes of
functions: For $ 0\le a<b\le\infty$,
\begin{enumerate}
\item$ \C[a,b)$: Continuous functions on $ [a,b)$.
\item$ \D[a,b)$: Right-continuous functions with finite left limits
defined on $ [a,b)$.
\item$ \D_{ l}[a,b)$: Left-continuous functions with finite right
limits defined on $ [a,b)$.
\end{enumerate}
$\D[0,1]$ is complete and separable under a metric $d_0(\cdot)$, which
is equivalent to the Skorohod metric $d_S(\cdot)$ (Billingsley~\cite{billingsley1968}, page~128),
but not under the uniform metric $ \|\cdot\|$. As we will see, the
limit processes that appear in our analysis below are always
continuous. We can
check that if $x$ is continuous (in fact, uniformly continuous) in
$[0,1]$, for $x_n \in\D[0,1]$, $\Vert x_n-x\Vert  \to0$ is equivalent to
$d_S(x_n,x) \to0$
and hence equivalent to $d_0(x_n,x) \to0$ as $n \to\infty$
(Billingsley~\cite{billingsley1968}, page~124). So we use
convergence in uniform metric, for our convenience
henceforth. For spaces of the form $ \D[a,b)$ or $ \D_{ l}[a,b)$, we
will consider the topology of locally uniform convergence. In some
cases we will also consider
product spaces of functions, and then the topology will be the product
topology. For example, $ \D^{ 2}_{ l}[1,\infty)$ will denote the class
of 2-dimensional functions on $ [1,\infty)$ which are left-continuous
with right limit. The classes of functions defined on the sets $ [a,b]$
or $ (a,b]$ will have the obvious notation.

\subsection{A useful lemma}
The following lemma will be used often in the proofs below. We use ``$
\Rightarrow$'' to denote weak convergence.

\begin{Lemma}\label{lemfnconvtosetconv}
Let $ Y_{ n}\in\D^{ 2}_{ l}(0,1]$ be a sequence of random functions
and assume the following hold:\vadjust{\goodbreak}
\begin{enumerate}[(ii)]
\item[(i)]$ Y_{ n}\Rightarrow Y$, where $ Y(t)$ has continuous paths with
probability 1.
\item[(ii)] There exists a partition $ 0=t^{ (0)}_{ n}<t_{ n}^{ (1)}<\cdots
<t_{ n}^{ (m_{ n})}=1$ such that $ Y_{ n}(t)$ is constant on the
interval $ (t_{ n}^{ (i)},t_{ n}^{ (i+1)}]$ for all $ 0\le i< m_{ n}$
with probability 1.
\end{enumerate}
Then for any $ 0<\varepsilon<1$,
%
%
\begin{equation}\label{eqfnconvtosetconvpart1}
\mathcal{Y}^{ \varepsilon}_{ n}:= \bigl\{ Y_{ n}\bigl(t_{ n}^{ (i)}\bigr)\dvt 0< i \le m_{
n},t_{ n}^{ (i)}\ge\varepsilon\bigr\}\quad\Rightarrow\quad\mathcal{Y}^{ \varepsilon
}:=\{ Y(t)\dvt \varepsilon\le t\le1 \} \qquad\mbox{in } \mathcal{F}.
\end{equation}
Furthermore, if $ \lim_{ t \downarrow0,n\to\infty} |Y_{ n}(t)|
=\infty$ with probability 1, then
%
%
\begin{equation}\label{eqfnconvtosetconvpart2}
\mathcal{Y}_{ n}:= \bigl\{ Y_{ n}\bigl(t_{ n}^{ (i)}\bigr)\dvt 0< i \le m_{ n} \bigr\}
\quad\Rightarrow\quad\mathcal{Y}:=\{ Y(t)\dvt 0< t\le1 \}\qquad \mbox{in }
\mathcal{F}.
\end{equation}
\end{Lemma}

\begin{pf}
Using Lemma~\ref{thmchoq} it suffices to show that
\[
\lim_{ n \rightarrow\infty}E[ h^{ \vee}(\mathcal{Y}_{ n}) ]
= E[ h^{ \vee}(\mathcal{Y}) ]
\]
for any continuous function $ h\dvtx\mathbb{R}^{ 2}\mapsto\mathbb{R}_{ +}$
with a compact support. So take any such function $ h$.
By the Skorohod representation theorem (Billingsley~\cite{billingsley1968}, Theorem~6.7, page~70), there exists a probability space $
(\Omega,\mathcal{G},P)$ and random elements
$ Y_{ n}^{ *}$ and $ Y^{ *}$ in $ \D^{ 2}_{ l}(0,1]$ such that
\[
Y_{ n}^{ *}\stackrel{ d}{= } Y_{ n}\quad \mbox{and}\quad Y^{
*}\stackrel{ d}{= } Y
\]
in the sense of finite dimensional distributions (f.d.d.) and
\[
Y_{ n}^{ *} \to Y^{* }, \qquad P\mbox{-a.s. in } \D^{ 2}_{ l}(0,1].
\]
Now observe that
\[
h^{ \vee}(\mathcal{Y}^{ \varepsilon})=\sup_{ x\in\mathcal{Y}^{ \varepsilon
}}h(x) \stackrel{ d}{= } \sup_{\varepsilon\le t\le1} h(Y^{ *}(t)
) \quad\mbox{and}\quad h^{ \vee}(\mathcal{Y}^{ \varepsilon}_{
n})= \sup_{ x\in\mathcal{Y}^{ \varepsilon}_{ n}}h(x) \stackrel{ d}{= }
\sup_{ \varepsilon\le t\le1} h(Y_{ n}^{ *}(t)).
\]
Since $ Y^{ *}(t)$ is continuous, we know that $ \sup_{ \varepsilon\le
t\le1} | Y^{ *}_{ n}(t)-Y^{ *}(t)|\to0$. Moreover, since $ h$ is
continuous with a compact support, we get $ h$ is uniformly continuous,
and hence
\begin{equation}\label{eqhYn*tohY}
\sup_{ \varepsilon\le t\le1} h(Y_{ n}^{ *}(t)) \to\sup_{
\varepsilon\le t\le1} h(Y^{ *}(t)), \qquad P\mbox{-a.s.}
\end{equation}
As $ h$ is bounded, applying the dominated convergence theorem, we get
\[
E[ h^{ \vee}(\mathcal{Y}^{ \varepsilon}_{ n}) ]= E \Bigl[ \sup_{
\varepsilon\le t\le1} h(Y_{ n}^{ *}(t)) \Bigr] \rightarrow E
\Bigl[ \sup_{ \varepsilon\le t\le1} h(Y^{ *}(t)) \Bigr] =E[ h^{
\vee}(\mathcal{Y}^{ \varepsilon}) ],
\]
and this proves \eqref{eqfnconvtosetconvpart1}.

Since $ h\dvtx \R^2 \mapsto\R_+$ has a bounded support, we can find $ M>0$
such that $ h(x)=0$ whenever $ |x|>M$. If $ \lim_{ t \downarrow0,n\to
\infty} |Y^{ *}_{ n}(t)| =\infty$ with probability 1, then almost
surely for any $ \omega\in\Omega$ we can find $ \delta>0$ and $
N\ge1$ such that $ |Y^{ *}_{ n}(t)(\omega)|>M$ for all $\delta\le
t\le1$ and $n\ge N$. This implies (using \eqref{eqhYn*tohY})
\[
\sup_{ 0<t\le1} h(Y_{ n}^{ *}(t)(\omega)) = \sup_{ \varepsilon\le t\le
1} h(Y_{ n}^{ *}(t)(\omega)) \to\sup_{ \varepsilon\le t\le1} h(Y^{
*}(t)(\omega))= \sup_{ 0< t\le1} h(Y^{ *}(t)(\omega)).\vadjust{\goodbreak}
\]
The remaining part of the proof of \eqref{eqfnconvtosetconvpart2} can
be completed using the same argument used to prove~\eqref
{eqfnconvtosetconvpart1}.
\end{pf}

\section{Limit results for the QQ plots} \label{secQQplot}
Convergence of empirical processes and quantile processes to
functionals of Gaussian processes, usually Brownian motion and Brownian
bridges, are quite well known; cf. Shorack and Wellner~\cite
{shorackwellner1986}. We
prove similar results for extreme order statistics. We use the weak
limit of tail empirical measure and deduce weak convergence of the
logarithmic version of the QQ plot of the extreme order statistics as a
random set.

The following was proved in Das and Resnick~\cite{dasresnick2008}:
\begin{Proposition}\label{propngenF}
Suppose $X_1,\ldots, X_n$ are i.i.d. with common distribution
$F$, and $ X_{(1)} \ge X_{(2)} \ge\cdots\ge X_{(n)}$
are the order statistics from this sample.
If $F$ is strictly increasing and continuous on its support, then
\[
\mathcal{T}_n : =\biggl\{\biggl(F^{\leftarrow}\biggl(\frac{i}{n+1}
\biggr), X_{(n-i+1)}\biggr) \dvt 1 \le i
\le n\biggr\}\quad \stackrel{P}{\to}\quad \mathcal{T} := \{ (x,x) \dvt x
\in \operatorname{support}(F)\}
\]
in $\F$.
\end{Proposition}
This proposition though is not enough if one is interested in creating
confidence bounds from the data.
For that purpose one would need weak convergence results which are
widely known in terms of
convergence of affine transformations of quantile processes to
appropriate Brownian Bridges for a known distribution $F$; see Shorack
and Wellner~\cite{shorackwellner1986}, Chapter~3, for further details.
In the following section, we concentrate on the case where $\bar{F}$ is
regularly varying with tail index $-1/\xi$ with $\xi>0$. The specific
form of $F$ is otherwise unknown.

\subsection{QQ plots for distributions with regularly varying tails}
Now assume that $X_1,X_2,\ldots,X_n$ are i.i.d. from a distribution
$F$. Suppose we want to check whether $F$ is heavy-tailed or not.
In the sense of testing a hypothesis,
our null hypothesis is that $\bar{F} \in RV_{-1/\xi}$ for some $\xi>0$.
Note that we really do not have any specific form for $F$. We define
the following sets:
\begin{eqnarray}
\mathcal{Q}_n & =& \biggl\{\biggl( -\log\frac jk, \log\frac
{X_{(j)}}{X_{(k)}} \biggr)\dvt 1\le j \le k \biggr\},\qquad k <n, \label
{eqnSnqq}\\
\mathcal{Q} & =& \{( x, \xi x )\dvt x \ge0 \}. \label{eqnTQQ}
\end{eqnarray}
The set $ \mathcal{Q}_{ n}$ is the logarithmic version of the QQ plot
for the first $ k$ order statistics from the sample $ X_{ 1},\ldots,X_{ n}$.

Das and Resnick~\cite{dasresnick2008} proved that under the null hypothesis,
$\mathcal{Q}_n \cinP\mathcal{Q}$ in $ \mathcal{F}$ as $k,n \to\infty$
with $k/n \to0$. We show below that a distributional convergence can
also be obtained in this case.
%
\begin{Assumption} \label{assmpQQ}
$F$ satisfies
\begin{equation}
\lim_{n\to\infty} \sqrt{k} \biggl( \frac nk \bar{F}
\bigl(b(n/k)y^{ -\xi}\bigr) - y\biggr)  = 0
\end{equation}
locally uniformly on $(0,\infty]$ as $k,n,n/k \to\infty$.
\end{Assumption}

\begin{Theorem}\label{thmdistQQRV}
Suppose $ X_{ 1},\ldots,X_{ n}$ are i.i.d. observations from a
distribution $ F$ satisfying $\bar{F} \in RV_{-1/\xi}$ with $\xi>0$
and Assumption~\ref{assmpQQ}. Then as $n,k,n/k\to\infty$
\begin{eqnarray}\label{eqQN}
&&\mathcal{QN}_{ n}:= \biggl\{ \biggl( -\log\frac{ j}{k },-\xi\log\frac
{ j}{k }+ \sqrt{k}\biggl( \log\frac{ X_{ (j)}}{ X_{ (k)}} +\xi\log\frac
{ j}{k } \biggr)\biggr)\dvt 1\le j\le k \biggr\}
\nonumber
\\[-8pt]
\\[-8pt]
\nonumber
&&\quad \Rightarrow\quad\mathcal{QN}:= \bigl\{ \bigl(-\log t,-\xi\log t +\xi t^{
-1}B(t) \bigr)\dvt 0<t\le1 \bigr\} \qquad\mbox{in } \mathcal{F},
\end{eqnarray}
where $B(t)$ is a Brownian Bridge on $[0,1]$ restricted to $(0,1]$.
\end{Theorem}

\begin{Remark}\label{remQQplot}
The set $ \mathcal{QN}_{ n}$ is a suitably normalized version of the
QQ plot which allows us to obtain a weak limit.
It is important to observe that the format in which we have expressed
the result is not standard in the literature as far as weak limits
of random variables or functions are concerned. Usual weak limit
results will only consider the normalized difference of the random variable
from its mean or its limit in probability. In our setting it is
imperative to state the result in the form which we have used. We look at
the plot as the probability limit perturbed by the normalized deviation
around it; that is, we shift the normalized differences so that
we can obtain the distribution of the deviation of the observed points
of the QQ plot from its mean position.
If we do not make this shift, the weak limit will always hover around
the $ y$-axis and will not give the deviation from the actual point in
the plot.
\end{Remark}

\begin{Remark}\label{remconditionR}
We have used Assumption~\ref{assmpQQ}
in order to prove a weak limit for the QQ plots. Without this
assumption we can show the convergence of tail empirical
measure with unknown centering $\frac{ n}{k } \bar F(b(n/k)y^{ -\xi})$
as in \eqref{eqtailconvthm91}, but we wish the centering to be $y$
here. To achieve this
\[
\lim_{n \to\infty}\sqrt{k}\biggl(\frac{ n}{k } \bar F\bigl(b(n/k)y^{
-\xi}\bigr)-y\biggr)
\]
should exist and have a finite value which we assume to be $ 0$ without
loss of any generality. The same theorem can be proved by replacing
Assumption~\ref{assmpQQ}
with the stronger condition of \textit{second order regular variation};
see de~Haan and Ferreira~\cite{dehaanferreira2006}, de~Haan and
Stadtmueller~\cite{dehaanstadtmueller1996},
de~Haan and Peng~\cite{dehaanpeng1998}.
Neither Assumption~\ref{assmpQQ} nor the second order RV condition is
easy to check in practice, albeit we resort to assuming them in order
to obtain distributional limits.

\end{Remark}

\begin{pf*}{Proof of Theorem~\protect\ref{thmdistQQRV}} The tail empirical measure defined as
%
%
\begin{equation}\label{eqtailempirical1}
\nu_{ n}(\cdot):= \frac{ 1}{ k}\sum_{ i=1}^{ n} \varepsilon_{ X_{
i}/b(n/k)}(\cdot)
\end{equation}
is a random element of $ M_{ +}(0,\infty]$ where $ \varepsilon_{ x}(\cdot
)$ puts unit mass at $ x$. By Theorem 4.1 (Resnick~\cite{resnickbook2007}, page~79), we get that
%
%
\begin{equation}\label{eqmeasureweak}
\nu_{ n} \Rightarrow\nu \qquad\mbox{in } M_{ +}(0,\infty],
\end{equation}
where $ \nu(y,\infty]=y^{ -1/\xi}, y>0.$
Furthermore, Theorem 9.1 in Resnick~\cite{resnickbook2007}, page~292,
gives us
%
%
\begin{equation}
\sqrt{k} \biggl( \nu_{ n} (y^{ -\xi},\infty]- \frac{ n}{k } \bar
F\bigl(b(n/k)y^{ -\xi}\bigr) \biggr)\Rightarrow W(y)\qquad \mbox{in }\D_{
l}(0,\infty], \label{eqtailconvthm91}
\end{equation}
where $W$ is a standard Brownian motion on $[0,\infty)$.
Since $ F$ satisfies Assumption~\ref{assmpQQ}, we obtain
%
%
\begin{equation}\label{eqtailemplimit2}
\sqrt{k} \bigl( \nu_{ n} (y^{ -\xi},\infty]- y\bigr)\Rightarrow W(y)\qquad \mbox{in }\D
_{ l}(0,\infty].
\end{equation}
We will use this to find the limiting distribution of
\[
\sqrt{k} \biggl( \log\frac{ X_{ (\lceil kt \rceil)}}{ X_{ (k)}} +\xi
\log t \biggr) = \sqrt{k} \log\biggl(\frac{X_{(\lceil kt\rceil
)}}{X_{(k)}} t^{\xi}\biggr),\qquad 0<t \le1,
\]
where for any $z\in\R$, denote by $\lceil z \rceil$, the largest
integer less than or equal to $z$. For $0 < t \le1$, let
\begin{eqnarray*}
\nu_n^{\la}(t) := \inf\{y\dvt \nu_n(y^{-\xi},\infty] \ge t\}
= \inf\Biggl\{y\dvt \sum_{ i=1}^{ n} \varepsilon_{ X_{ i}/b(n/k)} (y^{-\xi
},\infty]\ge kt\Biggr\}
= \biggl(\frac{X_{(\lceil kt\rceil)}}{b(n/k)}\biggr)^{-1/\xi}.
\end{eqnarray*}
Note that we can apply Vervaat's lemma (Resnick~\cite{resnickbook2007}, Proposition~3.3, page~59) to
\eqref{eqtailemplimit2} to get
%
%
\begin{equation}\label{eqvervaat}
\sqrt{k} \biggl( \biggl(\frac{ X_{(\lceil k t \rceil)}}{b(n/k)}
\biggr)^{ -1/\xi}- t \biggr) \Rightarrow W(t) \qquad\mbox{in } \D_{
l}(0,1] .
\end{equation}
Therefore, using the continuous map $f\dvtx\D_{l}(0,1] \to\D_{l}(0,1]$
with $f(x)(t) = x(t)/t$, we have
\begin{equation}\label{eqvervaatplus}
\sqrt{k} \biggl( \biggl(\frac{ X_{ (\lceil k t \rceil)}}{b(n/k)} t^{\xi}
\biggr)^{ -1/\xi}- 1 \biggr)  \Rightarrow\frac{W(t)}t\qquad \mbox{in } \D_{ l}(0,1].
\end{equation}
Also observe that
\begin{eqnarray}\label{eqlogtononlogforqq}
\sqrt{k} \log\biggl( \frac{ X_{ (\lceil k t \rceil)}}{b(n/k)} t^{\xi
}\biggr) & =& - \sqrt{k}\xi\log\biggl[1- \biggl(1-\biggl(\frac{ X_{
(\lceil k t \rceil)}}{b(n/k)} t^{\xi}\biggr)^{-1/\xi}\biggr)\biggr]
\nonumber
\\
&= &- \sqrt{k}\xi\biggl( \biggl(\frac{ X_{ (\lceil k t \rceil)}}{b(n/k)}
t^{\xi} \biggr)^{ -1/\xi} -1 \biggr) \\
&&{}+ \mathrm{o}_P\biggl(\sqrt{k}\xi\biggl(
\biggl(\frac{ X_{ (\lceil k t \rceil)}}{b(n/k)} t^{\xi} \biggr)^{ -1/\xi}-1
\biggr)\biggr).\nonumber
\end{eqnarray}
So from \eqref{eqvervaatplus} and \eqref
{eqlogtononlogforqq} it follows that
\begin{equation}
\sqrt{k} \log\biggl( \frac{ X_{ (\lceil k t \rceil)}}{b(n/k)} t^{\xi
}\biggr)  \Rightarrow- \xi\frac{W(t)}t \qquad\mbox{in } \D_{
l}(0,1].\label{eqvervaatinlogterms}
\end{equation}
We again use the continuous mapping theorem with $f\dvtx\D_{
l}(0,1] \to\D_{ l}(0,1]$, defined as $f(x)(t) = x(t) - x(1)$, to get
the following:
\begin{eqnarray} \label{eqmainqq1}
\sqrt{k} \log\biggl(\frac{ X_{ (\lceil k t \rceil)}}{X_{ (k)}} t^{\xi
}\biggr) & = &-\sqrt{k} \log\frac{ X_{ ( k )}}{b(n/k)} + \sqrt{k} \log
\frac{ X_{ (\lceil k t \rceil)}}{b(n/k)}t^{ \xi}
\nonumber
\\[-8pt]
\\[-8pt]
\nonumber
& \Rightarrow&\xi W(1) -\xi\frac{W(t)}t\qquad \mbox{in } \D_{
l}(0,1].
\end{eqnarray}
%
We know that $tW(1)- W(t) \eqd B (t)$ on $\D_{ l}[0,1]$, where ``$ \eqd
$'' denotes equality in distribution, and~$B$ is a Brownian Bridge on
$[0,1]$. Therefore, it is true on a restriction, and hence
\[
\sqrt{k} \biggl( \log\frac{ X_{ (\lceil kt \rceil)}}{ X_{ (k)}} +\xi
\log t \biggr) \weak\xi t^{ -1} B(t)\qquad \mbox{in } \D_{ l}(0,1].
\]
Furthermore, we also get
\begin{eqnarray}\label{eqfunctionformlimit}
&&{S}_{ n}(t)= \biggl( -\log\frac{ \lceil kt \rceil}{k },-\xi\log
\frac{ \lceil kt \rceil}{ k}+\sqrt{k} \biggl( \log\frac{ X_{ (\lceil kt
\rceil)}}{ X_{ (k)}} +\xi\log\frac{ \lceil kt \rceil}{k } \biggr)\biggr)
\nonumber
\\[-8pt]
\\[-8pt]
\nonumber
&&\quad \Rightarrow\quad {S}(t)= \biggl(-\log t,-\xi\log t +\xi\frac{B(t)}{t} \biggr)\qquad \mbox{in } \D
^{ 2}_{ l}(0,1],
\end{eqnarray}
using the converging-together lemma (Resnick~\cite{resnickbook2007}, Proposition~3.1, page~57) and the fact that
\[
\sqrt{k} \biggl( \log\frac{ \lceil kt \rceil}{ k} -\log t \biggr) \to0
\]
locally uniformly on $ (0,1]$. The weak convergence of the set $
\mathcal{QN}_{ n}$ follows from Lemma~\ref{lemfnconvtosetconv} once we
note that $ S_{ n}$ and $ S$ in \eqref{eqfunctionformlimit} satisfy
the conditions of Lemma~\ref{lemfnconvtosetconv}.
\end{pf*}


\section{Limit results for the ME Plots}\label{secMEplot}

\subsection{Empirical ME function for known distribution $F$}
Suppose $ X_{ 1},\ldots,X_{ n}$ is an i.i.d. sample from distribution $
F$. Yang~\cite{yang1978} studied the properties of the empirical ME function
$ \hat M(u)$ in \eqref{eqempiricalmeanexcess} as an estimator of $
M(u)$. They showed that $ \hat M(u)$ is uniformly strongly consistent for
$ M(u)$: for any $ 0<b<\infty$,
\[
P\Bigl[ \lim_{ n \rightarrow\infty} \sup_{ 0\le u\le b} | \hat
M(u)-M(u) |= 0 \Bigr]=1.
\]
Yang~\cite{yang1978} also proved a weak limit for $ \hat M(u)$: for any
$ 0<b<1$,
\[
\sqrt{n} \bigl( \hat M(F^{ \leftarrow}(t)) -M(F^{
\leftarrow}(t)) \bigr) \Rightarrow U(t),
\]
where $U(t) $ is a Gaussian process on $ [0,b]$ with covariance function
\[
\Gamma(s,t)= \frac{ (1-t)\sigma^{ 2}(t)-t\theta^{
2}(t)}{(1-s)(1-t)^{ 2} }\qquad \mbox{for all }0\le s\le t\le b
\]
with
\[
\sigma^{ 2}(t)=\operatorname{var} \bigl( XI_{ [t < F(X)\le1]} \bigr)\quad \mbox{and}\quad \theta(t)=E \bigl( XI_{ [t < F(X)\le1]} \bigr) .
\]
Although these properties are stated for the empirical ME function,
using Lemma~\ref{lemfnconvtosetconv}, it can be shown that the ME
plots also exhibit the same features when the distribution $ F$ is known.

\subsection{ME plot in the regularly varying case}
The behavior of $ \hat M(u)$ near the right end-point of $ F$ is not
explained in Yang~\cite{yang1978}. Here we study the asymptotic properties
of the ME plot when the explicit form of the distribution $ F$ is not known.
Ghosh and Resnick~\cite{ghoshresnick2009} proved the limit in
probability of a
suitably scaled version of the ME plot under the following null hypothesis:
\begin{Theorem}\label{thmpositivexi}
If $ X_{ 1},\ldots,X_{ n}$ are i.i.d. observations with distribution $
F$ satisfying $ \bar F\in RV_{ -1/\xi}$ with $ 0<\xi<1$, then in $
\mathcal{F}$,
%
%
\begin{equation}\label{eqlimpos}
\mathcal{M}_{ n}:=\frac{ 1}{X_{ (k)} } \bigl\{ \bigl(X_{ (i)},\hat
M\bigl(X_{ (i)}\bigr)\bigr) \dvt i=2,\ldots, k\bigr\} \quad\stackrel{ P}{\longrightarrow}\quad
\mathcal{M}:=\biggl\{ \biggl(t, \frac{ \xi}{1-\xi} t \biggr)\dvt t\ge1 \biggr\} .
\end{equation}
\end{Theorem}

In this paper we obtain the weak limit of the ME plot when the null
hypothesis that $ \bar F\in RV_{ -1/\xi}$ for some $ \xi>0$ holds.
The limit distribution depends on the value of $ \xi$. We get different
limits depending on whether $ \xi\le1/2$, $ 1/2<\xi<1$ or $ \xi\ge1$.

\subsubsection{\texorpdfstring{Case I: $0<\xi<1/2$}{Case I: 0<xi<1/2}}
In this case $ \operatorname{var} (X_{ 1}) < \infty$ exists and we obtain a Gaussian
limit for the suitably normalized ME plots. The following assumption is
essential. It is stronger than Assumption~\ref{assmpQQ} which was
required to obtain the weak limit of the QQ plot. As we discussed in
Remark~\ref{remconditionR}, it is quite difficult to check this
assumption in practice.
\begin{Assumption}\label{assmpME1}
$ F$ satisfies Assumption~\ref{assmpQQ}, and, moreover,
\[
\sqrt{k} \int_{ 1}^{ \infty}\biggl|\frac nk \bar{F} \bigl(b(n/k)y\bigr) -
y^{ -1/\xi} \biggr| \,\mathrm{d}y\to0
\]
as $n,k,n/k \to\infty$.
\end{Assumption}

\begin{Theorem}\label{thmweaklimitMEplot}
Suppose $ X_{ 1},\ldots,X_{ n}$ are i.i.d. observations from a
distribution $ F$ satisfying $ \bar F\in RV_{ -1/\xi}$ with $ 0<\xi
<1/2$ and Assumption~\ref{assmpME1} holds. Then for any $ 0<\varepsilon
<1$, as $ n,k,n/k\to\infty$,
\begin{eqnarray*}
&&\mathcal{MN}_{ n}:=\biggl\{\biggl( \biggl(\frac{ i}{k }\biggr)^{ -\xi},
\frac{ \xi}{1-\xi} \biggl( \frac{ i}{k } \biggr)^{ -\xi} \biggr) \\
&&\phantom{\mathcal{MN}_{ n}:=\biggl\{}{} + \sqrt{k}\biggl( \frac{ X_{ (i)}}{X_{
(k)} } -\biggl(\frac{ i}{k }\biggr)^{ -\xi}, \frac{\hat M(X_{ (i)})}{X_{
(k)}} - \frac{ \xi}{1-\xi} \biggl( \frac{ i}{k } \biggr)^{ -\xi} \biggr)\dvt i=\lceil\varepsilon k\rceil,\ldots,k\biggr\} \\
&&\quad \Rightarrow\quad \mathcal{MN}:= \biggl\{\biggl( t^{ -\xi}+\xi t^{ -(1+\xi
)}B(t),\frac{ \xi}{1-\xi}t^{ -\xi}+ \xi t^{ -1} \int_{ 0}^{ t} y^{
-(1+\xi)}B(y)\,\mathrm{d}y \biggr) , \\
&&\quad \phantom{\Rightarrow\quad \mathcal{MN}:= \biggl\{}
\varepsilon\le t\le1\biggr\} \qquad\mbox{in
} \mathcal{F},
\end{eqnarray*}
where $ B(t)$ is the standard Brownian bridge on $ [0,1]$ restricted to
$ (0,1]$.
\end{Theorem}

\begin{Remark}
Similar to Theorem~\ref{thmdistQQRV} we look at the ME plot as the
probability limit perturbed by the normalized deviation around it and
obtain a weak limit in Theorem~\ref{thmweaklimitMEplot}. The
assumption that $ \xi<1/2$ is essential. Note that
\[
\int_{ 0}^{ t} y^{ -(1+\xi)}W(y)\,\mathrm{d}y= \int_{ t^{ -\xi}}^{ \infty}W(u^{
-1/\xi})\,\mathrm{d}u=\int_{ t^{-\xi}}^{ \infty}\int_{ 0}^{ y^{ -1/\xi}}\,\mathrm{d}W(s)\,\mathrm{d}y=
\int_{ 0}^{ t}s^{ -\xi}\,\mathrm{d}W(s),
\]
and it is well known that the integral on the right-hand side exists if
and only if $ \int_{ 0}^{ t}s^{ -2\xi}\,\mathrm{d}s< \infty$, for which it is
necessary and sufficient to have $ \xi<1/2$; cf. {\O}ksendal
\cite{oksendal2003}, Lemma~3.1.5,
page~26. This means
\[
\int_{ 0}^{ t} y^{ -(1+\xi)}B(y)\,\mathrm{d}y \stackrel{ d}{= } \int_{ 0}^{ t} y^{
-(1+\xi)}W(y)\,\mathrm{d}y-W(1)\int_{ 0}^{ t}y^{ -\xi} \,\mathrm{d}y
\]
exists if and only if $ \xi<1/2$, and the same is true for the limit $
\mathcal{MN}$.
\end{Remark}

\begin{pf*}{Proof of Theorem \protect\ref{thmweaklimitMEplot}}
Consider a functional form of the ME plot,
%
%
\begin{equation}\label{eqMEplotfuncform}
S_{ n}(t)= \bigl( S_{ n}^{ (1)}(t),S_{ n}^{ (2)}(t) \bigr) :=\biggl(
\frac{ X_{ (\lceil k t \rceil)}}{X_{ (k)} }, \frac{\hat M(X_{(\lceil
kt\rceil) }) }{X_{ (k)}} \biggr) ,\qquad t\in(0,1],
\end{equation}
as random elements in $ \D^{ 2}_{ l}(0,1]$.
Following the proof of Theorem 3.2 in Ghosh and Resnick \cite
{ghoshresnick2009}, we
know that $S_{ n}(\cdot) \stackrel{ P}{ \longrightarrow} S(\cdot) $ in
$ \D^{ 2}_{ l}(0,1]$, where
\[
S(t):=\bigl(S^{ (1)}(t),S^{ (2)}(t)\bigr)= \biggl(t^{ -\xi}, \frac{ \xi}{1-\xi}
t^{ -\xi} \biggr),\qquad t\in(0,1].
\]
Applying Vervaat's lemma (Resnick~\cite{resnickbook2007}, Proposition~3.3, page~59) to \eqref{eqtailemplimit2}, we get
%
%
\begin{eqnarray}\label{eqvervaatmeplot}
&&\biggl(\sqrt{k} \biggl( \biggl(\frac{ X_{(\lceil k t \rceil)}}{b(n/k)}
\biggr)^{ -1/\xi}- t \biggr) ,\sqrt{k} \bigl( \nu_{ n}(t^{ -\xi},\infty
]-t \bigr) \biggr)
\nonumber
\\[-8pt]
\\[-8pt]
\nonumber
&&\quad\Rightarrow(-W(t),W(t))\qquad \mbox{in } \D^{ 2}_{
l}(0,\infty) .
\end{eqnarray}
Observe that
\[
S^{ (2)}_{ n}(t):= \frac{\hat M(X_{(\lceil kt \rceil)})}{X_{ (k)}} =
\frac{ k}{\lceil kt\rceil-1 }\int_{ X_{ (\lceil kt\rceil)}/b(n/k)}^{
\infty} \hat\nu_{ n}(x, \infty]\,\mathrm{d}x,
\]
where
%
%
\begin{equation}\label{eqnuhat}
\hat\nu_{ n}(\cdot):= \frac{ 1}{ k}\sum_{ i=1}^{ n} \varepsilon_{ X_{
i}/X_{ (k)}}(\cdot).
\end{equation}
Using \eqref{eqvervaatmeplot} and the converging-together lemma
(Resnick~\cite{resnickbook2007}, Proposition 3.1, page 57), we also have
\begin{eqnarray}\label{eqthreeconvmeplot}
&&\biggl(\sqrt{k} \biggl( \biggl(\frac{ X_{(\lceil k t \rceil)}}{b(n/k)}
\biggr)^{ -1/\xi}- t \biggr) ,\sqrt{k} \bigl( \nu_{ n}(t^{ -\xi},\infty
]-t \bigr), \frac{X_{(k)}}{b(n/k)} \biggr)
\nonumber
\\[-8pt]
\\[-8pt]
\nonumber
&&\quad\Rightarrow(-W(t),W(t),1) \qquad\mbox{in } \D
^{ 2}_{ l}(0,\infty)\times(0,\infty) .
\end{eqnarray}
Define a map $\hat{T}\dvtx\D^{ 2}_{ l}(0,\infty)\times(0,\infty) \to\D_{
l}(0,\infty)\times\D[1,\infty)$ as $\hat
{T}(f,g,x)(t,y)=(f(t), g(y^{-1/\xi}x)+y^{-1/\xi}f(1))$. We can check
that $\hat{T}$ is continuous at any $(f_0,g_0,x_0) \in\C^2(0,\infty
]\times(0,\infty)$ (see Resnick~\cite{resnickbook2007}, page 83). Hence,\
by the continuous mapping theorem on \eqref{eqthreeconvmeplot}, with
the map $\hat{T}$, we get that in $\D_{ l}(0,\infty)\times\D[1,\infty)$
\begin{eqnarray}\label{eqHk1}
J_n(t,y)&:=& \biggl(\sqrt{k} \biggl( \biggl(\frac{ X_{(\lceil k t \rceil
)}}{b(n/k)} \biggr)^{ -1/\xi}- t \biggr),\sqrt{k} \bigl( \hat\nu_{ n}
(y,\infty]- y^{-1/\xi} \bigr)\biggr) \nonumber\\
&=& \biggl(\sqrt{k} \biggl( \biggl(\frac{ X_{(\lceil k t \rceil)}}{b(n/k)}
\biggr)^{ -1/\xi}- t \biggr),
\nonumber
\\
&&\phantom{\biggl(} \sqrt{k} \biggl( \hat\nu_{ n} \biggl(y\frac
{X_{(k)}}{b(n/k)},\infty\biggr]- \biggl(y\frac{X_{(k)}}{b(n/k)}
\biggr)^{-1/\xi} \biggr)\\
&&\phantom{\biggl(}{}+ y^{-1/\xi}\sqrt{k}\biggl(\biggl(\frac
{X_{(k)}}{b(n/k)}\biggr)^{ -1/\xi}-1\biggr) \biggr)\nonumber \\
& \Rightarrow&\bigl(-W(t),W(y^{ -1/\xi})-y^{ -1/\xi}W(1)\bigr)
.\nonumber
\end{eqnarray}
{By an application of the functional delta method (van der Vaart and
Wellner~\cite{vandervaartwellner1996}, Theorem~3.9.4) to \eqref{eqHk1}, we obtain in $\D
^*:= \D_{ l}(0,1] \times\D[1,\infty)$,
%
%
\begin{eqnarray}\label{eqJk}
J^*_n(t,y)&:=&\biggl( \sqrt{k} \biggl( \frac{ X_{ (\lceil k t \rceil
)}}{b(n/k) } - t^{ -\xi} \biggr) ,\sqrt{k} \bigl( \hat\nu_{ n} (y,\infty
]- y^{-1/\xi} \bigr) \biggr)
\nonumber
\\[-8pt]
\\[-8pt]
\nonumber
&\Rightarrow&\bigl(\xi t^{ -(1+\xi)} W(t),
W(y^{ -1/\xi})-y^{ -1/\xi}W(1)\bigr).
\end{eqnarray}
The map $\phi\dvtx\D^* \to\D^* $ given by $\phi(f,g)=(f^{-\xi},g)$
is Hadamard differentiable at $(f(t),g(y))= (t,y^{-1/\xi})$,
tangentially to $\D_0:=\C(0,1]\times\C[1,\infty) \subset\D^*$, with
the right-hand side of \eqref{eqHk1}
being separable and an element of $\D_0$. Thus we can apply the
functional delta method (van der Vaart and Wellner~\cite{vandervaartwellner1996}, Theorem~3.9.4) to obtain \eqref{eqJk}.
}

Consequently, in $ \D_{ l}(0,1] \times\D[1,\infty)$,
\begin{eqnarray}\label{eqHk}
H_{ n}(t,y)
&:=&\biggl (\sqrt{k} \biggl( \frac{ X_{ (\lceil k t \rceil)}}{X_{ (k)} }
- t^{ -\xi} \biggr) ,\sqrt{k} \bigl( \hat\nu_{ n} (y,\infty]- y^{-1/\xi
} \bigr) \biggr)\nonumber\\
& = &\biggl(\frac{ b(n/k)}{X_{ (k)} }\sqrt{k} \biggl( \frac{ X_{ (\lceil
k t \rceil)}}{b(n/k) } - t^{ -\xi} \biggr) -\frac{ b(n/k)}{X_{ (k)}
}\sqrt{k} \biggl( \frac{ X_{ ( k )}}{b(n/k) }-1 \biggr)t^{ -\xi} ,\nonumber\\
&&\phantom{\biggl(}\sqrt
{k}\bigl ( \hat\nu_{ n} (y,\infty]- y^{-1/\xi} \bigr) \biggr)
\\
&\Rightarrow&\bigl(\xi t^{ -(1+\xi)} W(t)-\xi t^{ -\xi}W(1) , W(y^{
-1/\xi})-y^{ -1/\xi}W(1)\bigr)\nonumber \\
& \stackrel{ d}{= }& \bigl(\xi t^{ -(1+\xi)} B(t), B(y^{ -1/\xi})\bigr)=:
H(t,y).\nonumber
\end{eqnarray}
%
Define, for some $1\le K < \infty$, the maps $ T$ and $ T_{
K}$ from $ \D_{ l}(0,1] \times\D[1,\infty)$ to $ \D_{ l}(0,1] \times\D
[1,\infty)$ by
%
%
\begin{equation}\label{eqmapsTTk}
T(f,g)(t,y) = \biggl(f(t), \int_{ y}^{ \infty}g(x) \,\mathrm{d}x\biggr)\quad \mbox{and}\quad T_{ K}(f,g)(t,y) = \biggl( f(t) ,\int_{ y}^{ K\vee y}g(x) \,\mathrm{d}x \biggr).
\end{equation}
We understand $ \int_{ y}^{ \infty}g(x)\,\mathrm{d}x= \infty$ if $ g$ is not
integrable. Note that, in the Skorohod metric $d_S$,
we get $d_S(T_K(f_n,g_{ n}),T_K(f,g)) \le d_{ S}(f_{ n},f)+ K
d_S(g_n,g) \to0$ where $\{f_n, n\ge1\}, f \in\D_l(0,1]$ and $\{g_n,
n\ge1\}, g \in\D_l[1,\infty)$ with
$d_S(f_n,f) \to0$ and $d_S(g_n,g)\to0$ as $n\to\infty$. So $T_K$ is
a continuous mapping.
By \eqref{eqHk} and the continuity of the map $ T_{ K}$, we get that $
T_{ K}(H_{ n})\Rightarrow T_{ K}(H)$. We also claim that, for any $
\varepsilon>0$,
%
%
\begin{equation}\label{eqconvergingtogether}
\lim_{ K \rightarrow\infty} \limsup_{ n \rightarrow\infty} P
[\| T_{ K}(H_{ n}) -T(H_{ n})\| >\varepsilon]=0.
\end{equation}
Note that, for any $ \varepsilon>0$,
\begin{eqnarray*}
&&\lim_{ K \rightarrow\infty} \limsup_{ n \rightarrow\infty} P
[\| T_{ K}(H_{ n}) -T(H_{ n})\| >\varepsilon]\\
&&\quad\le\lim_{ K \rightarrow\infty} \limsup_{ n \rightarrow\infty}
P\biggl[\sqrt{k}\biggl|\int_{ K}^{ \infty}\bigl(\hat\nu_{ n}(x,\infty]
-x^{ -1/\xi}\bigr)\,\mathrm{d}x\biggr|>\varepsilon\biggr]\\
&&\quad \le\lim_{ K \rightarrow\infty} \limsup_{ n \rightarrow\infty}
P\biggl[\sqrt{k}\biggl|\int_{ K}^{ \infty}\biggl(\nu_{ n}\biggl(x\frac{ X_{
(k)}}{b(n/k) },\infty\biggr] -\biggl(x\frac{ X_{ (k)}}{b(n/k) }\biggr)^{
-1/\xi}\biggr)\,\mathrm{d}x\biggr|>\varepsilon/2 \biggr]\\
&&\qquad{} + \lim_{ K \rightarrow\infty} \limsup_{ n \rightarrow
\infty} P\biggl[\sqrt{k}\biggl|\int_{ K}^{ \infty}x^{ -1/\xi}\biggl(
\biggl(\frac{ X_{ (k)}}{b(n/k) }\biggr)^{ -1/\xi}-1\biggr)\,\mathrm{d}x\biggr|>\varepsilon/2
\biggr].
\end{eqnarray*}
Using \eqref{eqvervaat} and the assumption that $ \xi<1/2$, we get
\[
\lim_{ K \rightarrow\infty} \limsup_{ n \rightarrow\infty} P
\biggl[\sqrt{k}\biggl|\int_{ K}^{ \infty}x^{ -1/\xi}\biggl(\biggl(\frac{ X_{
(k)}}{b(n/k) }\biggr)^{ -1/\xi}-1\biggr)\,\mathrm{d}x\biggr|>\varepsilon/2 \biggr]=0.
\]
Using a change of variable, we obtain
\begin{eqnarray*}
&& \lim_{ K \rightarrow\infty} \limsup_{ n \rightarrow\infty}
P\biggl[\sqrt{k}\biggl|\int_{ K}^{ \infty}\biggl(\nu_{ n}\biggl(x\frac{ X_{
(k)}}{b(n/k) },\infty\biggr] -\biggl(x\frac{ X_{ (k)}}{b(n/k) }\biggr)^{
-1/\xi}\biggr)\,\mathrm{d}x\biggr|>\varepsilon/2 \biggr] \\
& &\quad= \lim_{ K \rightarrow\infty} \limsup_{ n \rightarrow\infty}
P\biggl[\sqrt{k}\biggl|\int_{ KX_{ (k)}/b(n/k)}^{ \infty}\bigl(\nu_{
n}(u,\infty] -u^{ -1/\xi}\bigr)\frac{ b(n/k)}{X_{ (k)}
}\,\mathrm{d}u\biggr|>\varepsilon/2 \biggr].
\end{eqnarray*}
Now fix any $ \eta>0$, and note that
\[
\lim_{ n \rightarrow\infty}P\biggl[ \biggl| \frac{ X_{ (k)}}{b(n/k) }-1
\biggr|>\eta\mbox{ or } \biggl| \frac{b(n/k) }{ X_{ (k)}}-1 \biggr|>\eta
\biggr] =0.
\]
Therefore,
\begin{eqnarray*}
&& \lim_{ K \rightarrow\infty} \limsup_{ n \rightarrow\infty}
P\biggl[\sqrt{k}\biggl|\int_{ KX_{ (k)}/b(n/k)}^{ \infty}\bigl(\nu_{
n}(u,\infty] -u^{ -1/\xi}\bigr)\frac{ b(n/k)}{X_{ (k)}
}\,\mathrm{d}u\biggr|>\varepsilon/2 \biggr] \\
& &\quad\le\lim_{ K \rightarrow\infty} \limsup_{ n \rightarrow\infty}
P\biggl[(1+\eta)\sqrt{k}\int_{ K(1-\eta)}^{ \infty}|\nu_{ n}
(u,\infty] -u^{ -1/\xi}|\,\mathrm{d}u>\varepsilon/2 \biggr] +\mathrm{o}(1).
\end{eqnarray*}
Now, since $ F$ satisfies Assumption~\ref{assmpME1}, it suffices to
show that
%
%
\begin{equation}\label{eqconvtogether}
\lim_{ K \rightarrow\infty} \limsup_{ n \rightarrow\infty} P
\biggl[\sqrt{k}\int_{ K(1-\eta)}^{ \infty}\biggl|\nu_{ n}(x,\infty] - \frac{
n}{k }\bar F\bigl( b(n/k)x \bigr) \biggr|\,\mathrm{d}x> \frac{ \varepsilon}{2(1-\eta) }
\biggr]= 0.
\end{equation}
This can be easily proved using the arguments in the proof of
Proposition 9.1 in (Resnick~\cite{resnickbook2007}, page 296). Observe
that, by using the triangle and Chebyshev inequalities,
\begin{eqnarray*}
&& P\biggl[\sqrt{k}\int_{ K(1-\eta)}^{ \infty}\biggl|\nu_{ n}(x,\infty] -
\frac{ n}{k }\bar F\bigl( b(n/k)x \bigr) \biggr|\,\mathrm{d}x>\frac{ \varepsilon
}{2(1-\eta) } \biggr] \\
&&\quad\le P\biggl[\frac{1}{\sqrt{k}}\int_{ K(1-\eta)}^{ \infty}\sum_{ i=1}^{
n}\bigl|\varepsilon_{ X_{ i}/b(n/k)}(x,\infty] - \bar F\bigl( b(n/k)x \bigr)
\bigr|\,\mathrm{d}x>\frac{ \varepsilon}{2(1-\eta) } \biggr]\\
& &\quad\le\biggl( \frac{ \varepsilon}{2(1-\eta) } \biggr) ^{- 2}\int_{
K(1-\eta)}^{ \infty}\frac{ n}{k }\operatorname{var} \bigl[\varepsilon_{ X_{
i}/b(n/k)}(x,\infty] \bigr] \,\mathrm{d}x\\
&&\quad\le\biggl( \frac{ \varepsilon}{2(1-\eta) } \biggr) ^{ -2} \int_{
K(1-\eta)}^{ \infty}\frac{ n}{k } \bar F \bigl( b(n/k)x \bigr)\,\mathrm{d}x \stackrel
{n\to\infty}{\to} \biggl( \frac{ \varepsilon}{2(1-\eta) } \biggr) ^{
-2}\int_{ K(1-\eta)}^{ \infty}x^{ -1/\xi}\,\mathrm{d}x.
\end{eqnarray*}
The last limit follows from Karamata's Theorem; cf. Resnick~\cite{resnickbook2007}, page~25.
Since $ \xi<1$, the integral in the last
expression is finite and therefore \eqref{eqconvtogether}, and hence
\eqref{eqconvergingtogether}, holds. From Theorem 3.5 in Resnick
\cite{resnickbook2007}, page~56, we get $ T(H_{ n})\Rightarrow T(H)=(\xi
t^{ -(1+\xi)} B(t),\int_{ y}^{ \infty} B(x^{ -1/\xi})\,\mathrm{d}x)$ in $ \D_{
l}(0,1] \times\D[1,\infty)$.

Now consider the random element $ Y_{ n}$ in the space $ \D^{ 2}_{
l}(0,1]\times\D[1, \infty)$,
%
\[
Y_{ n}(t,y):=\biggl( \frac{ X_{ (\lceil k t \rceil)}}{X_{ (k)} }, T(H_{
n})(t,y) \biggr).
\]
By another application of the converging-together lemma, it is easy to
check that $ Y_{ n}\Rightarrow Y $, where
\[
Y(t,y)=\biggl(t^{ -\xi}, \xi t^{ -(1+\xi)} B(t), \int_{ y}^{ \infty}
B(x^{ -1/\xi})\,\mathrm{d}x\biggr) .
\]
The map $ \tilde T\dvtx\D^{ 2}_{ l}(0,1]\times\D[1, \infty)\to\D^{ 2}_{
l}(0,1] $ defined by
\[
\tilde T\bigl(\bigl(f^{ (1)},f^{ (2)}\bigr),g\bigr)(t)=\bigl(f^{ (2)}(t),g\bigl(f^{
(1)}(t)\bigr)\bigr)\qquad \mbox{for all } 0<t\le1
\]
is continuous at $ (f,g)\in\C^{ 2}(0,1]\times\C[1,\infty)$. Therefore
\begin{eqnarray*}
\tilde T (Y_{ n})(t) &=& \biggl( \sqrt{k} \biggl( \frac{ X_{ (\lceil k t
\rceil)}}{X_{ (k)} } - t^{ -\xi} \biggr), \sqrt{k} \int_{ X_{ (\lceil
kt\rceil)}/X_{ (k)}}^{ \infty} \bigl( \hat\nu_{ n}(y, \infty]-y^{ -1/\xi
}\bigr)\,\mathrm{d}y \biggr)\\
&\Rightarrow&\biggl(\xi t^{ -(1+\xi)} B(t),\int_{ t^{ -\xi}}^{ \infty}
B(y^{ -1/\xi})\,\mathrm{d}y \biggr)\qquad \mbox{in } \D^{ 2}_{ l}(0,1].
\end{eqnarray*}
This implies
\[
\sqrt{k} \bigl( S_{ n}(t)-S(t) \bigr) \Rightarrow\biggl(\xi t^{ -(1+\xi)}
B(t) , \frac{ 1}{ t} \int_{ t^{ -\xi}}^{ \infty} B(y^{ -1/\xi})\,\mathrm{d}y\biggr)
\qquad\mbox{in } \D^{ 2}_{ l}(0,1].
\]
It is then easy to check that
\begin{eqnarray*}
&&\biggl(\xi t^{ -(1+\xi)} B(t) , \frac{ 1}{ t} \int_{ t^{ -\xi}}^{ \infty
} B(y^{ -1/\xi})\,\mathrm{d}y\biggr)
\\
&&\quad\stackrel{ d}{= } \biggl(\xi t^{ -(1+\xi)} B(t)
, \frac{ \xi}{ t} \int_{ 0}^{ t} y^{ -(1+\xi)}B(y)\,\mathrm{d}y\biggr)\qquad \mbox{in }\D^{ 2}_{ l}(0,1].
\end{eqnarray*}
%
Also observe that
\begin{eqnarray}\label{eqMEtildeS}
&&\tilde S_{ n}(t):= \biggl(\biggl( \frac{ \lceil kt \rceil}{ k} \biggr)
^{ -\xi}, \frac{ \xi}{1-\xi}\biggl ( \frac{ \lceil kt \rceil}{ k}
\biggr)^{ -\xi}\biggr)\\
&&\phantom{\tilde S_{ n}(t):=}{}
+\sqrt{k} \biggl( \frac{ X_{ (\lceil k t
\rceil)}}{X_{ (k)} }-\biggl( \frac{ \lceil kt \rceil}{ k} \biggr) ^{ -\xi
}, \frac{\hat M(X_{(\lceil kt\rceil) }) }{X_{ (k)}}- \frac{ \xi}{1-\xi
} \biggl( \frac{ \lceil kt \rceil}{ k} \biggr)^{ -\xi} \biggr)\\
&&\quad \Rightarrow\quad\tilde S(t):=\biggl(t^{ -\xi}, \frac{ \xi}{1-\xi}t^{
-\xi} \biggr)
\nonumber
\\[-8pt]
\\[-8pt]
\nonumber
&&\quad\phantom{\Rightarrow\quad\tilde S(t):= }{}+\biggl(\xi t^{ -(1+\xi)} B(t) , \frac{ \xi}{ t} \int_{ 0}^{
t} y^{ -(1+\xi)}B(y)\,\mathrm{d}y\biggr) \qquad\mbox{in }\D^{ 2}_{ l}(0,1],
\end{eqnarray}
since
\[
\sqrt{k}\biggl(\biggl( \frac{ \lceil kt \rceil}{ k}\biggr)^{ -\xi}-t^{ -\xi
} \biggr) \to0 \qquad\mbox{as } k\to\infty
\]
locally uniformly on $ (0,1]$. The proof the theorem is completed by
applying Lemma~\ref{lemfnconvtosetconv} to $ \tilde S_{ n}$ and $
\tilde S$.
\end{pf*}

\subsection{\texorpdfstring{Case II: $1/2<\xi<1$}{Case II: 1/2<xi<1}}

When $ 1/2<\xi<1$, the distribution $ F$ admits a finite mean but not a
finite variance. The ME function, however, exists, and we know the
limit in probability of the scaled ME plot from Theorem~\ref{thmpositivexi}.

\begin{Assumption}\label{assmpME2}
$ F$ satisfies Assumption~\ref{assmpQQ}, and, moreover,
\[
\frac{ 1}{b(n) } \biggl(\frac{k b(n/k)}{1-\xi}u^{1 -\xi}-C_{ ku,n}
\biggr)\to0
\]
for every $ 0<u<1$. For any $ l<n$
%
%
\begin{equation}\label{eqCkn}
C_{ l,n}:= n \int_{ 0}^{ l/n}F^{ \leftarrow}(1-u)\,\mathrm{d}u.
\end{equation}
\end{Assumption}

\begin{Theorem}\label{thmxigehalf}
Suppose $ X_{ 1},\ldots,X_{ n}$ are i.i.d. observations from a
distribution $ F$ satisfying $ \bar F\in RV_{ -1/\xi}$ with $ 1/2<\xi
<1$ and Assumption~\ref{assmpME2}. Then for any $ 0<\varepsilon<1$
\begin{eqnarray*}
&&\mathcal{MN}_{ n}:= \biggl\{\biggl(\biggl (\frac{ i}{k }\biggr)^{ -\xi},
\frac{ \xi}{1-\xi} \biggl( \frac{ i}{k } \biggr)^{ -\xi} \biggr) \\
 &&\phantom{\mathcal{MN}_{ n}:= \biggl\{}{}+
\biggl(\sqrt{k}\biggl( \frac{ X_{ (i)}}{X_{ (k)} } -\biggl(\frac{ i}{k }
\biggr)^{ -\xi}\biggr), \frac{ kb(n/k)}{b(n) } \biggl(\frac{\hat M(X_{
(i)})}{X_{ (k)}} - \frac{ \xi}{1-\xi} \biggl( \frac{ i}{k } \biggr)^{ -\xi
} \biggr) \biggr) \dvt\\
&&\phantom{\mathcal{MN}_{ n}:= \biggl\{}{} i=\lceil\varepsilon k \rceil,\ldots,k\biggr\} \\
&&\quad\Rightarrow\quad \mathcal{MN}:= \biggl\{\biggl( t^{ -\xi}+\xi t^{ -(1+\xi
)}B(t),\frac{ \xi}{1-\xi}t^{ -\xi}+ t^{ -1} S_{ 1/\xi} \biggr) ,
\varepsilon\le t\le1\biggr\}\qquad \mbox{in } \mathcal{F},
\end{eqnarray*}
where $ B(t)$ is the standard Brownian bridge on $ [0,1]$ restricted to
$ (0,1]$ and $ S_{ 1/\xi}$ is a stable random variable independent of
$ B(t)$ with characteristic function
%
%
\begin{equation}\label{eqstablegehalf}
E[ \mathrm{e}^{\mathrm{i}tS_{ 1/\xi}} ] =
\exp\biggl\{ -\frac{ 1}{1-\xi}\Gamma\biggl( 2- \frac{ 1}{\xi}
\biggr)\cos\frac{ \uppi}{2\xi}| t |^{ 1/\xi}\biggl[ 1-\mathrm{i}\sgn(t) \tan\frac{ \uppi
}{2\xi} \biggr] \biggr\}.
\end{equation}

\end{Theorem}

\begin{Remark}
An interesting point to note here is that the two coordinates of the
weak limit $ \mathcal{MN}$ are independent. The empirical ME function
depends on the sum of the order statistics $ X_{ (1)},\ldots, X_{
(k)}$. When $ 1/2<\xi<1$, this sum is dominated by a very few high
order statistics, and it turns out that the contribution of $ X_{ (k)}$
to the suitably normalized $\hat M(X_{ (k)}) $ vanishes in the limit.
The proof below formalizes this idea.

This feature is in stark contrast to what happens in the case $0< \xi
<1/2$. In that case all the top $ k$ order statistics have some
contribution to $ \hat M(X_{ (k)})$ in the limit. Hence the two
coordinates in the limit are obtained from the same Gaussian process
and are definitely not independent.
\end{Remark}

\begin{Remark}\label{remxieqhalf}
Unfortunately, we are unable to obtain a proper weak limit of the ME
plot in the case when $ \xi=1/2$. In this case it is known that the
weak limit of the suitably normalized sum of the first $ k$ order
statistics is Gaussian; cf. Cs\"orgo, Haeusler and Mason~\cite
{csorgohaeuslermason1991}. So this
would be similar to what happens when $0< \xi<1/2$, but the problem is
that the integral $ \int_{ 0}^{ t}y^{ -2}\,\mathrm{d}B(y)$ does not exist. It is
possible to redefine the ME plot in a different way, by leaving out a
few of the top order statistics and obtaining a limit in that case, but
we did not pursue that direction.
\end{Remark}

\begin{pf*}{Proof of Theorem \protect\ref{thmxigehalf}}
From Theorem 3 in Cs\"orgo, Horv\'ath and Mason \cite
{csorgohorvathmason1986} we know that if $
l=l_{ n}\to\infty$ with $ l_{ n}/n\to0$, then
%
%
\begin{equation}\label{eqcsorgostable}
\frac{ 1}{b(n) } \Biggl( \sum_{ i=1}^{ l}X_{ (i)}-C_{ l,n} \Biggr)
\Rightarrow S_{ 1/\xi}.
\end{equation}
Observe that, by Karamata's theorem (Resnick~\cite{resnickbook2007}, Theorem~2.1, page~25),
\[
C_{ l,n}=n\int_{ n/l}^{ \infty}b(s)/s^{ 2}\,\mathrm{d}s \sim n\frac{
(n/l)b(n/l)}{(n/l)^{ 2}(1-\xi) }= \frac{ lb(n/l)}{1-\xi}.
\]

Choose $ l=l_{ n}$ such that $ l/k\to0$ as $n\to\infty$. Fix any $
0<u<1$. Then
\begin{eqnarray}\label{eqgehalfpart1}
V_{ n}(t)& =&\bigl(V^{ (1)}_{ n}(t),V^{ (2)}_{ n}\bigr) : =\Biggl( \sqrt{k} \biggl(
\frac{ X_{ (\lceil k t \rceil)}}{X_{ (k)} } - t^{ -\xi} \biggr) ,\frac{
1}{b(n) } \Biggl( \sum_{ i=1}^{ l}X_{ (i)}-C_{ l,n} \Biggr)
\Biggr)
\nonumber
\\[-8pt]
\\[-8pt]
\nonumber
& \Rightarrow&\bigl(\xi t^{ -(1+\xi)} B(t),S_{ 1/\xi}\bigr)
\qquad\mbox{in } \D_{ l}[u,1] \times\mathbb{R},
\end{eqnarray}
where $ B(t)$ and $ S_{ 1/\xi}$ are as described in the statement of
the theorem. The convergence of the coordinates $ V^{ (1)}_{ n}(t)$ and
$ V^{ (2)}_{ n}$ of $ V_{ n}(t)$ follows from \eqref{eqHk} and \eqref
{eqcsorgostable}. The asymptotic independence of $ V^{ (1)}_{ n}(t)$
and $ V^{ (2)}_{ n}$ is a consequence of Theorem D in Cs\"orgo and
Mason \cite
{csorgomason1985} or Satz 4 in Rossberg~\cite{rossberg1967}.
Using \eqref{eqtailemplimit2} we get that
\[
\sqrt{k} \Biggl( \frac{ 1}{kb(n/k) } \sum_{ i= \lceil ku \rceil+1}^{
\lceil kt \rceil}X_{ (i)}- \frac{ 1}{1-\xi} ( t^{ 1-\xi}-u^{ 1-\xi
} ) \Biggr) \Rightarrow\int_{ u}^{ t}W(y)\,\mathrm{d}y \qquad\mbox{in } \D
_{ l}[u,1],
\]
and since $ \xi>1/2$, $ kb(n/k)/(b(n)\sqrt{k})\to0$, which implies
%
%
\begin{eqnarray}\label{eqhalfmidpart}
U^{ (2)}_{ n}(t)&:=& \frac{ kb(n/k)}{b(n) }\Biggl( \frac{ 1}{kb(n/k) }
\sum_{ i= \lceil ku \rceil+1}^{ \lceil kt \rceil}X_{ (i)}- \frac{
1}{1-\xi} ( t^{ 1-\xi}-u^{ 1-\xi} ) \Biggr)
\nonumber
\\[-8pt]
\\[-8pt]
\nonumber
 &\rightarrow&\mathbf
{0}\qquad \mbox{in } \D_{ l}[u,1],
\end{eqnarray}
where $ \mathbf{0}\in\D_{ l}[u,1]$ denotes the identically zero function.
Furthermore, using Theorem 2 in Cs\"orgo, Horv\'ath and Mason \cite
{csorgohorvathmason1986}, we get
\[
\frac{ 1}{\sqrt{k}b(n/k) } \Biggl( \sum^{ \lceil ku \rceil}_{i= l+1}X_{
(i)} -(C_{ ku,n}-C_{ l,n}) \Biggr)\Rightarrow N(0,1)
\]
and hence
%
%
\begin{equation}\label{equppermidpart}
U^{ (3)}_{ n}:=\frac{ 1}{b(n) } \Biggl( \sum^{ \lceil ku \rceil}_{i=
l+1}X_{ (i)} -(C_{ ku,n}-C_{ l,n}) \Biggr) \to0.
\end{equation}
Combining \eqref{eqgehalfpart1}, \eqref{eqhalfmidpart} and \eqref{equppermidpart} and the converging-together lemma
(Resnick~\cite{resnickbook2007},
Proposition 3.1, page 57), we get an important
building block of this proof,
%
%
\begin{eqnarray}\label{eqimportantblock}
U_{ n}(t)&:=& \bigl( V^{ (1)}_{ n}(t),U^{ (2)}_{ n}(t),U^{ (3)}_{ n},V^{
(2)}_{ n}\bigr)
\nonumber
\\[-8pt]
\\[-8pt]
\nonumber
&\Rightarrow&\bigl(\xi t^{ -(1+\xi)} W(t),\mathbf
{0},0,S_{ 1/\xi}\bigr)\qquad \mbox{in }\D^{ 2}_{ l}[u,1]\times\mathbb
{R}^{ 2}.
\end{eqnarray}

Next we consider
\begin{eqnarray*}
Z_{ n}(t)&= &\bigl( Z_{ n}^{ (1)}(t),Z_{ n}^{ (2)}(t) \bigr) := \biggl(
\sqrt{k} \biggl( \frac{ X_{ (\lceil k t \rceil)}}{X_{ (k)} } - t^{ -\xi}
\biggr) , \frac{ kb(n/k)}{b(n) } \biggl(\frac{\hat M(X_{ (\lceil kt
\rceil)})}{X_{ (k)}} - \frac{ \xi}{1-\xi} t^{ -\xi} \biggr) \biggr)
\\
&\in&\D_{ l}^{ 2}[u,1]
\end{eqnarray*}
and focus on the second coordinate $ Z_{ n}^{ (2)}(t)$.
\begin{eqnarray*}
Z_{ n}^{ (2)}(t) & =&\frac{ kb(n/k)}{b(n) } \biggl(\frac{\hat M(X_{
(\lceil kt \rceil)})}{X_{ (k)}} - \frac{ \xi}{1-\xi} t^{ -\xi}
\biggr)\\
& =& \frac{ kb(n/k)}{b(n) } \biggl(\frac{\hat M(X_{ (\lceil kt \rceil
)})}{b(n/k)} - \frac{ \xi}{1-\xi} t^{ -\xi} \biggr)+\mathrm{o}_{ P}(1)\\
& = &\frac{ kb(n/k)}{b(n) } \Biggl( \frac{ 1}{(\lceil kt \rceil-1)b(n/k)
}\sum_{ i=1}^{ \lceil kt \rceil-1}X_{ (i)}- \frac{ X_{ (\lceil kt \rceil
)}}{b(n/k) }- \frac{ \xi}{1-\xi} t^{ -\xi} \Biggr) +\mathrm{o}_{ P}(1)\\
& =& \frac{ kb(n/k)}{b(n) } \Biggl( \frac{ 1}{(\lceil kt \rceil-1)b(n/k)
}\sum_{ i=1}^{ \lceil kt \rceil-1}X_{ (i)}- \frac{ 1}{1-\xi} t^{ -\xi}
\Biggr) +\mathrm{o}_{ P}(1)\\
& =& \frac{ kb(n/k)}{tb(n) } \Biggl(\frac{ 1}{kb(n/k) } \sum_{ i=1}^{
\lceil kt \rceil-1}X_{ (i)}-\frac{ 1}{1-\xi} t^{ 1-\xi} \Biggr) +\mathrm{o}_{
P}(1)\\
& =& \frac{ 1}{t }U^{ (2)}_{ n}(t)+ \frac{ 1}{t } U^{ (3)}_{ n}+ \frac{
1}{t }V_{ n}^{ (2)}+\mathrm{o}_{ P}(1),
\end{eqnarray*}
where the last equality holds because of Assumption~\ref{assmpME2}.
Therefore, we get
\[
Z_{ n}(t)\Rightarrow\bigl(\xi t^{ -(1+\xi)} B(t),t^{ -1}S_{ 1/\xi}
\bigr) \qquad\mbox{in } \D^{ 2}_{ l}[u,1].
\]
Since the above limit holds for every $ 0<u<1$, it holds in $ \D_{
l}(0,1]$ as well. The proof is completed using Lemma~\ref{lemfnconvtosetconv}.
\end{pf*}

\subsection{\texorpdfstring{Case III: $\xi\ge1$}{Case III: xi>=1}}

In this case, the distribution $ F$ need not have a finite mean, and
the ME function may not be defined. It definitely does not exists if $
\xi>1$. Still the empirical ME plot can have a limit.

\begin{Theorem}\label{thmxiposgtr1}
Suppose $ X_{ 1},\ldots,X_{ n}$ are i.i.d. observations with
distribution $ F$ satisfying $ \bar F\in RV_{ -1/\xi}$ and Assumption
\ref{assmpME1}.
\begin{enumerate}
\item If $ \xi>1$ and $n,k,n/k \to\infty$, then
\begin{eqnarray*}
&&\mathcal{MN}_{ n}  := \biggl\{ \biggl(\biggl(\frac{ i}{k }\biggr)^{ -\xi}+\sqrt
{k}\biggl( \frac{ X_{ (i)}}{X_{ (k)} } -\biggl(\frac{ i}{k }\biggr)^{ -\xi
}\biggr),\frac{\hat M(X_{ (i)})}{b(n)/k}\biggr) \dvt i=2,\ldots, k \biggr\}
\\
&&\quad \Rightarrow\quad\mathcal{MN}:=\bigl\{ \bigl(\xi t^{ -(1+\xi)} B(t), tS_{
1/\xi} \bigr)\dvt t\ge1 \bigr\}
\end{eqnarray*}
in $ \mathcal{F}$, where $ S_{ 1/\xi}$ is the positive stable random
variable with index $ 1/\xi$ which satisfies, for $ t\in\mathbb{R}$,
\[
E[ \mathrm{e}^{\mathrm{i}tS_{ 1/\xi}} ] =
\exp\biggl\{ -\Gamma\biggl( 1- \frac{ 1}{\xi} \biggr)\cos\frac{ \uppi
}{2\xi}| t |^{ 1/\xi}\biggl[ 1-i\sgn(t) \tan\frac{ \uppi}{2\xi} \biggr]
\biggr\},
\]
and $ B(t)$ is a Brownian bridge independent of $ S_{ 1/\xi}$.
\item If $ \xi=1$, and $ k$ satisfies $n,k,n/k \to\infty$, and $
kb(n/k)/b(n)\to1$, then
\begin{eqnarray*}\label{eqlimposeql1}
&&\mathcal{MN}_{ n}:= \biggl\{ \biggl(\biggl(\frac{ i}{k }\biggr)^{ -\xi}+\sqrt
{k}\biggl( \frac{ X_{ (i)}}{X_{ (k)} } -\biggl(\frac{ i}{k }\biggr)^{ -\xi
}\biggr),\frac{\hat M(X_{ (i)})}{b(n/k)}- \frac{k
C_{ k,n}^{ *}}{ib(n) }\biggr) \dvt i=2,\ldots, k
\biggr\} \\
&&\quad \Rightarrow\quad\mathcal{MN}:=\bigl\{ t\bigl( t^{ -1} B(t), S_{ 1}-1-\log
t \bigr):t\ge1 \bigr\}
\end{eqnarray*}
in $ \mathcal{F}$, where
\[
C_{ k,n}^{ *} = n\int_{ 1/n}^{ k/n}F^{ \leftarrow}(1-u)\,\mathrm{d}u ,
\]
$ S_{ 1}$ is a positively skewed stable random variable satisfying
\[
E[ \mathrm{e}^{\mathrm{i}tS_{ 1}} ] =\exp\biggl\{ \mathrm{i} t\int_{ 0}^{ \infty} \biggl(
\frac{ \sin x}{ x^{ 2}}- \frac{ 1}{x(1+x) } \biggr)\,\mathrm{d}x - | t | \biggl[ \frac
{\uppi}{2}+\mathrm{i}\sgn(t) \log| t | \biggr] \biggr\},
\]
and $ B(t)$ is a Brownian bridge independent of $ S_{ 1}$.
\end{enumerate}

\end{Theorem}

\begin{pf}
The theorem is proved in the same fashion as the previous ones. First
we prove the weak limit in the functional form of the ME plot, and then
we infer the weak limit of the plot as a random set. Define
\[
S_{ n}(t)= \cases{
\displaystyle\biggl( \sqrt{k} \biggl( \frac{ X_{ (\lceil k t \rceil
)}}{X_{ (k)} } - t^{ -\xi} \biggr), \frac{ \hat M (X_{ (\lceil kt \rceil
)})}{ b(n)/k} \biggr) \vspace*{2pt}\cr
\quad\mbox{in part (i)}\vspace*{2pt}\cr
\displaystyle\biggl( \sqrt{k}
\biggl( \frac{ X_{ (\lceil k t \rceil)}}{X_{ (k)} } - t^{ -\xi} \biggr), \frac
{\hat M(X_{ (\lceil kt \rceil)})}{b(n/k)}- \frac{k
C_{ k,n}^{ *}}{\lceil kt \rceil b(n) }\biggr) \vspace*{2pt}\cr
\quad \mbox{in part
(ii)}}\qquad
\mbox{ for all }0<t\le1.
\]
We have already proved the weak limit of $ S_{ n}^{ (1)}(t)$ and the
weak limit of $ S_{ n}^{ (2)}(t)$ is proved in Theorem 3.4 in Ghosh and
Resnick \cite
{ghoshresnick2009}. The rest of the proof is completed using Lemma
\ref{lemfnconvtosetconv}.
\end{pf}

\section{Confidence bounds for the plots}\label{secconfidence}
In Sections~\ref{secQQplot} and~\ref{secMEplot}, we have obtained
weak convergence limits for the QQ and ME plots in the Fell topology.
Since the limit set in each case is a closed random set, we can compute
from the results in Sections~\ref{secQQplot} and~\ref{secMEplot}, the
probability that the random limit set
is contained in a fixed set in~$\R^2$.
This leads to creating asymptotic $ 100(1-\alpha) \%$ confidence
bounds around the plots, given any $0<\alpha<1$.
The methodology for creating confidence bounds around the plots is
explained in
details for QQ plots, and the same idea follows for ME plots.

\subsection{QQ plots}
Under the usual assumptions of Section~\ref{secQQplot}, the QQ plot,
$\mathcal{Q}_n$, as defined in \eqref{eqnSnqq}, consists of $k=k(n)<n$
points in $\R^2$. We know that $\mathcal{Q}_n\cinP\mathcal{Q}$,
where $\mathcal{Q}$ is a straight line.
From Theorem~\ref{thmdistQQRV}, we also know that $\mathcal{QN}_n$,
which is an affine transformation of $\mathcal{Q}_n$, converges weakly
to a random set $\mathcal{QN}$
centered around $\mathcal{Q}$ in $\mathcal{F}$. For fixed $0<\alpha
<1$, we intend to create a confidence bound around $\mathcal{Q}_n$
which will contain $\mathcal{Q}$
with probability $1-\alpha$ under the null hyothesis.

The limit distribution for QQ plots obtained in
Theorem~\ref{thmdistQQRV} is a linear transformation of $\{
t^{-1}{B(t)}\dvt 0< t \le1\}$ where $B$
is a Brownian bridge on $[0,1]$. So the limit explodes as $t$ comes
close to $0$, and thus we create confidence bounds under an $\varepsilon$
truncation to avoid this. Define
\begin{eqnarray}
\mathcal{Q}_n^{\varepsilon} & :=& \biggl\{\biggl( -\log\frac jk, \log\frac
{X_{(j)}}{X_{(k)}} \biggr)\dvt 1\le j \le k \mbox{ and } \frac jk\ge\varepsilon
\biggr\}, \qquad k <n\label{eqtruncQQplot} \\
\mathcal{Q}^{\varepsilon} & :=& \{( -\log t, -\xi\log t )\dvt
\varepsilon\le t \le1 \}=\biggl\{(x,\xi x)\dvt 0 \le x \le
\log\frac1{\varepsilon}\biggr\}\label{eqtruncQQlimit}.
\end{eqnarray}
Now with similar truncations defined as above, it follows from \eqref
{eqfunctionformlimit} that $\tilde{S}_n^{\varepsilon} \weak\tilde
{S}^{\varepsilon}$ in $\D^{ 2}_l[\varepsilon,1]$. This means $ \mathcal{QN}_{
n}^{ \varepsilon}\Rightarrow\mathcal{QN}^{ \varepsilon}$ in $ \mathcal{F}$
where $ \mathcal{QN}^{ \varepsilon}$ and $ \mathcal{QN}_{ n}^{ \varepsilon
}$ are the truncated versions of $ \mathcal{QN}$ and $ \mathcal{QN}_{
n}$, respectively, defined in \eqref{eqQN}.
Suppose we can calculate $c_{\alpha/2,\varepsilon}$ such that $P(\sup
_{\varepsilon\le t \le1} \frac{|B(t)|}{t} \le c_{\alpha/2,\varepsilon
} ) = 1-\alpha$.
Then a conservative $100(1-\alpha)\%$ confidence bound around $\mathcal
{Q}_n^{\varepsilon}$ is given by
\begin{equation}\label{eqconbandQQplot}
\mathcal{CQ}_n^{\varepsilon}= \mathcal{Q}_n^{\varepsilon} + \biggl\{
(0,y)\dvt y\in\xi\biggl(-\frac{c_{\alpha/2,\varepsilon}}{\sqrt{k}},\frac
{c_{\alpha/2,\varepsilon}}{\sqrt{k}}\biggr) \biggr\}.
\end{equation}

It is easy to see that
\[
P [ \mathcal{QN}^{ \varepsilon} \subset\mathcal{CQ}^{ \varepsilon}
] \ge1-\alpha,
\]
where
\[
\mathcal{CQ}^{ \varepsilon}:=\mathcal{Q}^{\varepsilon} + \biggl\{(0,y)\dvt y\in
\xi\biggl(-\frac{c_{\alpha/2,\varepsilon}}{\sqrt{k}},\frac{c_{\alpha
/2,\varepsilon}}{\sqrt{k}}\biggr) \biggr\}.
\]
An equivalent statement in a different notation is
\[
P [ \{ \rho(x,\mathcal{QN}^{ \varepsilon})=0 ,\forall x\in
\mathcal{CQ}^{ \varepsilon} \}] \ge1-\alpha,
\]
where, for any $ x\in\mathbb{R}^{ 2}$ and $ F\in\mathcal{F}$,
\[
\rho(x,F):=\inf\{ |x-y|\dvt y\in F \} .
\]
From (Molchanov~\cite{molchanov2005}, Theorems B.6 and B.13, pages 400--401),
we know that if $ F_{ n}\to F$ in~$ \mathcal{F}$, then for any compact
set $ K\subset\mathbb{R}^{ 2}$
\[
\sup_{ x\in K} | \rho(x,F_{ n}) -\rho(x,F)| \to0.
\]
Since $ \mathcal{QN}_{ n}^{ \varepsilon}\Rightarrow\mathcal{QN}^{
\varepsilon}$ in $ \mathcal{F}$, we get
\[
\lim_{ n \rightarrow\infty} P [ \{ \rho(x,\mathcal{QN}_{
n}^{ \varepsilon})=0 ,\forall x\in\mathcal{CQ}^{ \varepsilon} \}
] \ge1-\alpha.
\]
Hence, $ \mathcal{CQ}_{ n}^{ \varepsilon}$ in \eqref{eqconbandQQplot}
is an asymptotic $ 100(1-\alpha)\%$ confidence bound for $ \mathcal
{Q}^{ \varepsilon}$.

We calculate $P(\sup_{\varepsilon\le t \le1} \frac{|B(t)|}{t} \le
M)$ next in order to complete the construction.
Since $W(t) := (1+t) B (\frac{t}{t+1}), 0 \le t < \infty$ is
a Brownian motion on $[0,\infty)$, we can check that
\begin{equation}
\sup_{\varepsilon\le t \le1} \frac{|B(t)|}{t}  = \sup_{t
\ge\delta} \frac{|W(t)|}{t}, \label{eqBtWt}
\end{equation}
where $\delta=\frac{\varepsilon}{1-\varepsilon}$. In the following theorem we
compute the boundary-crossing probability for $\sup_{t \ge\delta
} \frac{|W(t)|}{t}$.
\begin{Proposition} \label{propquantileofWtbyt}
Suppose $W$ is a standard Brownian motion on $[0,\infty)$. Then for all
$\delta>0$ and $M>0$,
\begin{equation}
 P\biggl(\sup_{t \ge\delta} \frac{|W(t)|}{t} > M \biggr) = 4\sum
_{k=1}^{\infty}\bigl[\Phi\bigl((4k+1)M\sqrt{\delta}\bigr)-\Phi\bigl((4k-1)M\sqrt{\delta
}\bigr)\bigr], \label{eqquantileofWtbyt}
\end{equation}
where $\Phi(\cdot)$ denotes the c.d.f. of a standard normal distribution.
\end{Proposition}

\begin{pf}
We begin by observing that
\begin{eqnarray}\label{eqreadyfordoob}
&& P\biggl(\sup_{t \ge\delta} \frac{|W(t)|}{t} \le M \biggr)
\nonumber
\\[-8pt]
\\[-8pt]
\nonumber
&&\quad= \int
_{s=-M\delta}^{M\delta} P\bigl(-Mt-(M\delta+s) \le W(t) \le Mt
+ (M\delta-s), \forall t \ge0 \bigr) f_{W(\delta)}(s) \,\mathrm{d}s,
\end{eqnarray}
where $f_{W(\delta)}$ denotes the density of $W(\delta)$. The
right-hand side is obtained by conditioning on $W(\delta)=s$ and using
the fact that $\{W(t)-W(\delta)\dvt t\ge{\delta}\}$ is independent of
$W(\delta)$ and $\{W(\delta)\dvt t\ge\delta\} \eqd\{W(t)\dvt t\ge0\}.$
Now the above boundary (non-)crossing probablity of the Brownian
motion can be calculated using Doob~\cite{doob1949}, equation (4.3), as
\begin{eqnarray*}
&&P\biggl(\sup_{t \ge\delta} \frac{|W(t)|}{t} > M \biggr) \\
&&\quad = 1 -
\int_{s=-M\delta}^{M\delta}\Biggl[1 - \sum_{k=1}^{\infty}
(\mathrm{e}^{-2A_k} + \mathrm{e}^{-2B_k}-\mathrm{e}^{-2C_k}-\mathrm{e}^{-2D_k})\Biggr] f_{W(\delta
)}(s) \,\mathrm{d}s,
\end{eqnarray*}
where
\begin{eqnarray*}
A_k &=& [(2k-1)M]^2\delta- (2k-1)Ms, \qquad B_k  = [(2k-1)M]^2\delta+
(2k-1)Ms,\\
C_k &= &4k^2M^2\delta-2kMs, \qquad D_k  = 4k^2M^2\delta+2kMs.
\end{eqnarray*}
Since ${W(\delta)}\sim N(0,\delta)$, for any $a,b\in\R$, we have
\begin{equation}
\int_{-a}^a \mathrm{e}^{bs} f_{W(\delta)}(s)\,\mathrm{d}s = \mathrm{e}^{b^2\delta/2}
\biggl[\Phi\biggl(\frac{a-b\delta}{\sqrt{\delta}}\biggr) - \Phi\biggl(\frac
{-a-b\delta}{\sqrt{\delta}}\biggr) \biggr]. \label{eqimpeq}
\end{equation}
Now using \eqref{eqimpeq}, we can compute, for each $k\ge1$,
\begin{eqnarray*}
\int_{-M\delta}^{M\delta} (\mathrm{e}^{-2A_k}+\mathrm{e}^{-2B_k}
)f_{W(\delta)}(s)\,\mathrm{d}s & = &2 \bigl[\Phi\bigl((4k-1)M\sqrt{\delta}\bigr)-\Phi\bigl((4k-3)M\sqrt
{\delta}\bigr)\bigr],\\
\int_{-M\delta}^{M\delta} (\mathrm{e}^{-2C_k}+\mathrm{e}^{-2D_k}
)f_{W(\delta)}(s)\,\mathrm{d}s & =& 2 \bigl[\Phi\bigl((4k+1)M\sqrt{\delta}\bigr)-\Phi\bigl((4k-1)M\sqrt
{\delta}\bigr)\bigr].
\end{eqnarray*}
Therefore we get
\begin{eqnarray*}
&& P\biggl(\sup_{t \ge\delta} \frac{|W(t)|}{t} > M \biggr)\\
&&\quad = 1 -
\int_{s=-M\delta}^{M\delta}\Biggl[1 - \sum_{k=1}^{\infty}
(\mathrm{e}^{-2A_k} + \mathrm{e}^{-2B_k}-\mathrm{e}^{-2C_k}-\mathrm{e}^{-2D_k})\Biggr] f_Z(s) \,\mathrm{d}s\\
&&\quad = 4\sum_{k=1}^{\infty}\bigl[\Phi\bigl((4k-1)M\sqrt{\delta}\bigr)-\Phi
\bigl((4k-3)M\sqrt{\delta}\bigr)\bigr].
\end{eqnarray*}
\upqed\end{pf}

%
\begin{Remark}
Observe that the confidence bound in \eqref{eqconbandQQplot} depends
on the value of $ \xi$. While obtaining the width of the band, we
replace $ \xi$ by its Hill estimate (Resnick~\cite{resnickbook2007}, page 74). We could use any consistent estimator of $ \xi$
and the choice of the estimator does not seem to be important as far as
the simulation study is concerned. It is well known that estimating the
parameter $ \xi$ can often be extremely tricky, see ``Hill--Horror
plots'' in (Resnick~\cite{resnickbook2007}, page 87). But as far as
obtaining confidence bounds is concerned, we can get past that by using
a conservative estimate of~$ \xi$, that is, a value which we strongly
believe is not less than the true value of $ \xi$.
\end{Remark}

\begin{Remark}
It is clear that the probability calculated in Proposition \ref
{propquantileofWtbyt} is very close to 1 if $M\sqrt{\delta}$ is small.
We can approximate the infinite sum in \eqref{eqquantileofWtbyt} by a
finite sum whose limit depends on our choice of $M$ and $\delta$. We use
Proposition~\ref{propquantileofWtbyt} to create confidence bands for
the QQ plots in the examples in Section~\ref{secdata}. Simulation
suggests that considering the first 15 terms of the infinite sum is
enough to give us approximations
correct up to six decimal places.\vspace*{-2pt}
\end{Remark}

\begin{Remark}
It is possible to join the subsequent points in $\mathcal{Q}_n$ to make
a continuous curve $\mathcal{Q}^*_n \in\mathcal{F}$, and we can check
that $ \mathcal{Q}^{ *}_{ n}$ will converge to the same limit as that
of $ \mathcal{Q}_{ n}$ as $ n\to\infty$.
We mention this result here without proof, which can be completed
following Theorem~\ref{thmdistQQRV}.\vspace*{-2pt}
\end{Remark}


\subsection{ME plots}\vspace*{-2pt}
In Section~\ref{secMEplot} we obtained weak
limits for the ME plots, under the assumption that $\bar{F}\in
RV_{-1/\xi}$ with $\xi>0$, where $F$ denotes the underlying
distribution. We observed three separate limits in three different cases.

For the case $0<\xi< 1/2$, where $F$ has a finite second moment, we
obtain a limit in terms of functionals
of Brownian bridges (see Theorem~\ref{thmweaklimitMEplot}). In order
to convert this result to obtain confidence bounds, we need to compute
boundary-crossing probabilities for these functionals.
Analytical solution for such probabilities are available for linear
boundaries (Doob~\cite{doob1949}) and piecewise linear boundaries
(P\"otzelberger and Wang~\cite{potzelbergerwang2001}) in case of
Brownian motion on
$[0,\infty)$.
Probabilities for nonlinear boundaries, which happens to be our case,
are usually approximated using results for piecewise linear boundaries.
Instead of such approximations, we resort to Monte Carlo simulation to
find appropriate confidence bounds; see Section~\ref{secdata}.

For the case $1/2<\xi<1$, $F$ has a finite first moment, but its
second moment does not exist. The limit distribution for the affinely
transformed ME plot consists of a functional of a Brownian bridge
in the first component and a Stable distribution in the second
component. The feature here is that the normalization required to get
the limit depends on $b(n)$ and $b(n/k)$, which in turn depends on the
distribution function $ F$ and is hence unknown.
These can be estimated in practice with $X_{(1)}$ and $X_{(k)}$
respectively. Although, to justify such a procedure we would need to
know the joint behavior $(X_{(1)},X_{(k)}, \sum_{i=1}^k X_{(i)})$ when
$k,n$ and $n/k \to\infty$.
Results in Darling~\cite{darling1952}, Chow and Teugels~\cite{chowteugels1978}, Resnick~\cite{resnick1986}, Section 4,
are quite useful here. Using Theorem 5.3 in
Darling~\cite{darling1952} we can show that under the assumptions of Theorem
\ref{thmxigehalf},
\begin{eqnarray}\label{eqthm6mod}
&&\widetilde{\mathcal{MN}}_{ n}:= \biggl\{\biggl( \biggl(\frac{ i}{k }
\biggr)^{ -\xi}, \frac{ \xi}{1-\xi} \biggl( \frac{ i}{k } \biggr)^{ -\xi}
\biggr)\nonumber \\[-2pt]
& &\phantom{\widetilde{\mathcal{MN}}_{ n}:= \biggl\{}{}+ \biggl(\sqrt{k}\biggl( \frac{ X_{ (i)}}{X_{ (k)} } -\biggl(\frac{
i}{k }\biggr)^{ -\xi}\biggr), \frac{ kX_{ (k)}}{X_{ (1)} } \biggl(\frac
{\hat M(X_{ (i)})}{X_{ (k)}} - \frac{ \xi}{1-\xi} \biggl( \frac{ i}{k }
\biggr)^{ -\xi} \biggr) \biggr) \dvt
\nonumber
\\[-10pt]
\\[-10pt]
\nonumber
&&\phantom{\widetilde{\mathcal{MN}}_{ n}:= \biggl\{}{} i=2,\ldots,k\biggr\} \\
&&\quad\Rightarrow\quad \widetilde{\mathcal{MN}}:= \biggl\{\biggl( t^{ -\xi}+\xi
t^{ -(1+\xi)}B(t),\frac{ \xi}{1-\xi}t^{ -\xi}+ t^{ -1} \tilde S_{
1/\xi} \biggr) , 0<t\le1\biggr\}\qquad \mbox{in } \mathcal
{F},\nonumber
\end{eqnarray}
where $ \tilde S_{ 1/\xi}$ is independent of $B(t) $, and its
characteristic function is of the form
%
%
\begin{equation}\label{eqcharfunUV}
E[\mathrm{e}^{\mathrm{i}\lambda\tilde S_{ 1/\xi}}]= \mathrm{e}^{\mathrm{i}\lambda} \biggl( 1+ \frac{
i\lambda}{1-\xi}-\frac{ 1}{\xi}\int_{ 0}^{ 1}(\mathrm{e}^{\mathrm{i}t\lambda
}-1-it\lambda) t^{ -1-1/\xi}\,\mathrm{d}t\biggr) ^{ -1}.
\end{equation}
We again resort to Monte Carlo simulation to obtain confidence bounds
for the ME plots.

For $\xi\ge1$, $F$ need not have a finite mean, and the ME plot does
not have a non-trivial non-random limit. We obtain weak limits here in
Theorem~\ref{thmxiposgtr1}.
Clearly, calculating confidence bounds is not sensible here.
\subsubsection{Confidence bound for ME plots} We need to truncate the
ME plot near infinity in this case, since the weak limits we obtain
(Theorems~\ref{thmweaklimitMEplot}
and~\ref{thmxigehalf}) blow up there (relates to $t$ near $0$ in
the limit $\mathcal{MN}_n$).
According to~\eqref{eqlimpos}, $\mathcal{M}_n$ denotes the ME plot
for a sample of size $n$ (with $k<n$ top order statistics under
consideration). Define its truncated version
%
%
\begin{eqnarray}\label{eqlimposepsilon}
\mathcal{M}^{ \varepsilon}_{ n}&:=&\frac{ 1}{X_{ (k)} } \bigl\{ \bigl(X_{
(i)},\hat M\bigl(X_{ (i)}\bigr)\bigr) \dvt i=\lceil k\varepsilon\rceil,\ldots, k\bigr\}
\quad\mbox{and}
\nonumber
\\[-8pt]
\\[-8pt]
\nonumber
\mathcal{M}^{ \varepsilon}&:=&\biggl\{ \biggl(t,
\frac{ \xi}{1-\xi} t \biggr)\dvt \varepsilon\le t \le1 \biggr\} .
\end{eqnarray}
Then $ \mathcal{M}_{ n}^{ \varepsilon}\stackrel{ P}{\to} \mathcal{M}^{
\varepsilon} $.

{If $ 0<\xi<1/2$, then using Theorem~\ref{thmweaklimitMEplot} we
can give the $(1-\alpha)100 \%$ confidence band for $\mathcal
{M}^{\varepsilon}$ as
\begin{equation}\label{eqconbandMEplot}
\mathcal{CM}_n^{\varepsilon} := \mathcal{M}_n^{\varepsilon} + \biggl\{
(x,y)\dvt x\in\biggl(-\frac{c_{\alpha_1/2,\varepsilon}}{\sqrt{k}},\frac
{c_{\alpha_1/2,\varepsilon}}{\sqrt{k}}\biggr), y\in\biggl( -\frac{d_{\alpha
_2/2,\varepsilon}}{\sqrt{k}},\frac{d_{\alpha_2/2,\varepsilon}}{\sqrt{k}}
\biggr) \biggr\},
\end{equation}
where $\alpha_{1},\alpha_{2}>0$ is such that $\alpha=\alpha
_{1}+\alpha_{2} $ and
\begin{eqnarray}\label{eqcalpheps}
c_{\alpha,\varepsilon} &=& \mbox{$(1-\alpha)$th quantile of } \sup
_{\varepsilon\le t \le1}\xi t^{-(1+\xi)}B(t),
\nonumber
\\[-8pt]
\\[-8pt]
\nonumber
d_{\alpha,\varepsilon}& =& \mbox{$(1-\alpha)$th quantile of } \sup
_{\varepsilon\le t \le1}\xi t^{-1} \int_0^t y^{-(1+\xi
)}B(y)\,\mathrm{d}y.
\end{eqnarray}
$\mathcal{CM}_n^{\varepsilon} $ in
\eqref{eqconbandMEplot} provides an asymptotic confidence bound
around $\mathcal{M}^{\varepsilon}$ with $P(\mathcal{M}^{\varepsilon} \subset
\mathcal{CM}_n^{\varepsilon}) \ge(1-\alpha)$ for large $n$.}

If $ 1/2<\xi<1$, then we use Theorem~\ref{thmxigehalf} and its
modified form in \eqref{eqthm6mod} to give the $(1-\alpha)100\% $
confidence band for $ \mathcal{M}^{ \varepsilon}$ as
\begin{eqnarray}\label{eqconbandMEplotxigehalf}
\mathcal{CM}_n^{\varepsilon} &=&\biggl\{ \biggl( \frac{ X_{ (\lceil kt \rceil
)}}{ X_{ (k)}}, \frac{ \hat M(X_{ (\lceil kt \rceil)})}{ X_{ (k)} }
\biggr)
\nonumber
\\[-8pt]
\\[-8pt]
\nonumber
&&\phantom{\biggl\{}{}+ \biggl(-\frac{c_{\alpha_1/2,\varepsilon}}{\sqrt{k}},\frac
{c_{\alpha_1/2,\varepsilon}}{\sqrt{k}}\biggr)\times\biggl( \frac{X_{ (1)}
d_{1-\alpha_2/2}}{\lceil kt \rceil X_{ (k)}},\frac{X_{ (1)}d_{\alpha
_2/2}}{ \lceil kt \rceil X_{ (k)}} \biggr)\dvt \varepsilon\le t \le1
\biggr\},
\end{eqnarray}
where
\[
d_{\alpha} = \mbox{$(1-\alpha)$th quantile of } \tilde S_{ 1/\xi}
\mbox{ defined in \eqref{eqcharfunUV}}.
\]
Here $0<\alpha_1,\alpha_{ 2}<1$ are chosen such that $(1-\alpha
)=(1-\alpha_1)(1-\alpha_{ 2})$. Since the random components in the
first and second components in the limit (Theorem~\ref{thmxigehalf})
are independent, this gives us the right confidence interval so that
$P(\mathcal{M}^{\varepsilon} \subset\mathcal{CM}_n^{\varepsilon} ) \ge
1-\alpha$. The above quantiles are calculated by Monte Carlo methods
for the simulation we report in Section~\ref{secsimulme}.

\begin{Remark}
Throughout the literature of extreme value theory, the top $ k$ order
statistics where $ k =k_{ n}\to\infty$ and $ k/n\to0$ as $n \to\infty
$ is considered for inference.
The idea is that as the size of data increases we concentrate more on
the extreme right-hand tail of the underlying distribution. In practice though,
given a data set of fixed size $n$, albeit large, it is difficult to
decide on which value of $ k$ to choose. The popular solution is to try
out different values of $ k$;
see Embrechts, Kl\"{u}ppelberg and Mikosch~\cite{embrechtskluppelbergmikosch1997}, Chapter 6,
and Resnick~\cite{resnickbook2008}, Chapter 4, for further discussions on this issue.

In order to obtain confidence bounds for QQ plots and ME plots, along
with the problem of choosing $ k$, we also have to choose $ \varepsilon$.
The choice of $ \varepsilon$ should be such that, for the purpose of
drawing any inference, we leave out the region where data is sparse. In
practice, we have to try out different values of $ \varepsilon$ depending
on the size of the data and the choice of $ k$.
\end{Remark}
\begin{Remark}
An important point to note here is that we are suggesting to use the
weak limit of the QQ plot to obtain the confidence band. In practice,
even if we have a large data set, it will always be finite. A natural
question that arises here is what is the rate of convergence in these
cases. We do not have the answer at the moment, but all the simulation
studies that we have done strongly suggest that this method works well.
\end{Remark}
%
\section{QQ plot and ME plot in practice}\label{secdata}
\subsection{Simulation}
We do a simulation study using the software R to check how well this
method of obtaining confidence bounds for the QQ plot and the ME plot works.
\subsubsection{QQ plots}\label{secsimulqq}

\begin{figure}

\includegraphics{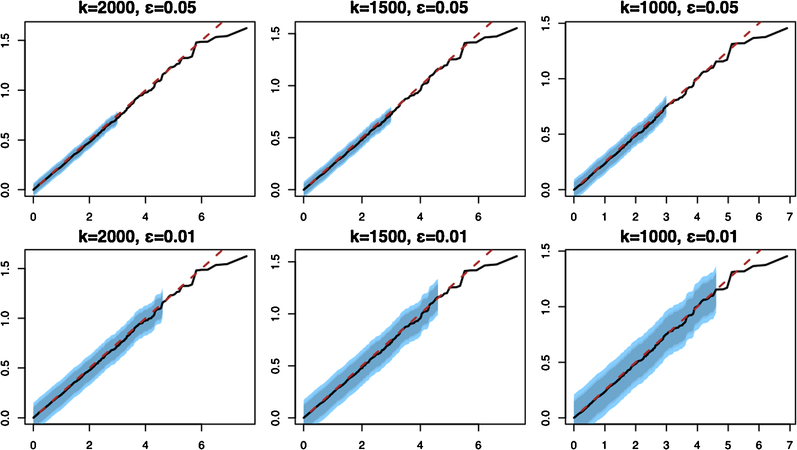}

\caption{QQ plot for 50,000 i.i.d. Pareto random variables with $ \xi
=0.25$.}\label{figQQpareto}
\end{figure}

We begin with a simple exercise for Pareto distribution with $ \xi
=0.25$ ($ \bar F(x)=x^{ -4}, x\ge1$). We simulate a sample of size $
n= 50\mbox{,}000$ from this distribution and look at the QQ plot for extremes
as defined in \eqref{eqnSnqq}; see Figure~\ref{figQQpareto}. The
black line denotes the plot $ \mathcal{Q}_{ n}$, and the brown dotted
line denotes the true line $ \mathcal{Q}$. We know that $ \mathcal{Q}_{
n}$ converges to $ \mathcal{Q}$, and, as we see in the plot, the two
lines are close, except for the top-right corner of the plots, which
correspond to the very large order statistics. We choose three
different values for $ k$: 2000, 1500 and 1000, which are large in
absolute terms, but small compared to the sample size $ n$.

%

Following the discussion in Section~\ref{secconfidence}, we know that
the variance of the limiting distribution blows up as we move towards
the extreme order statistics (towards the top-right corner) in the
plot. So while obtaining a confidence bound, we truncate at $ \lfloor
\varepsilon k\rfloor$th order statistic for $ \varepsilon=0.05 \mbox{ and
}0.01$. The confidence bounds are obtained for the six cases. The three
shades of the colored bands signify the $ 99\%, 95\%$ and the $ 90\%$
confidence bands for the plot. As is evident in Figure \ref
{figQQpareto}, the true line lies within the bound in all the cases.
It is also notable that the width of the confidence band increases as $
k$ and $ \varepsilon$ decrease.

%

\begin{figure}

\includegraphics{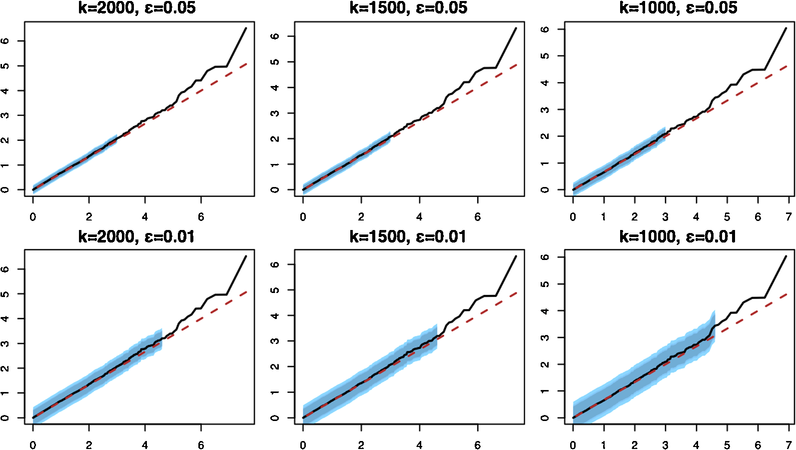}

\caption{QQ plot for 50,000 i.i.d. right-skewed stable random variables
with $ \xi=2/3$.}\label{figQQstable}\vspace*{-3pt}
\end{figure}

Next we do a similar study for a right-skewed stable distribution with
$ \xi=2/3$ $(\alpha=1.5)$ and mean 0. We use the same values of $ n,k$
and $ \varepsilon$. The result is given in Figure~\ref{figQQstable}.
Here also we see that the method works well, and the confidence band
contains the true line in all the six cases.

\begin{figure}

\includegraphics{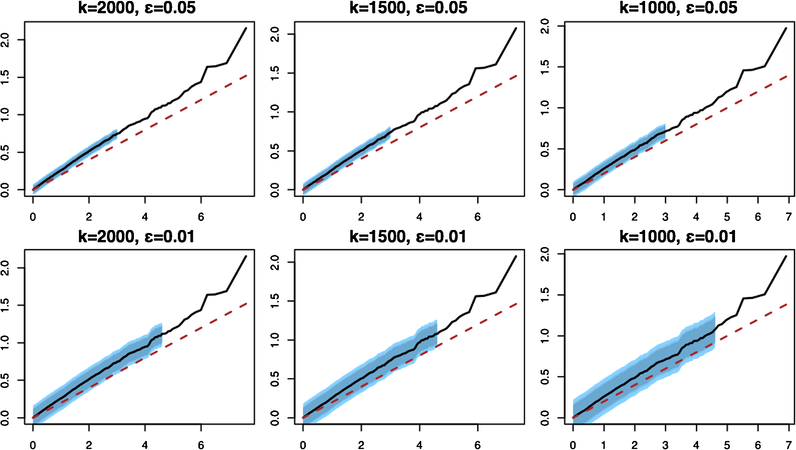}

\caption{QQ plot for 50,000 i.i.d. random variables with the
distribution described in \protect\eqref{eqnonstd} ($\xi=0.2$).}\label
{fignonstandard}
\end{figure}

We also try a non-standard distribution for which $ \bar F ^{
-1}(x)=x^{ -1/5}(1-10^{ -1}\ln x)$, $0<x\le1$. This means that $ \bar
F\in RV_{ -4}$, and therefore $ \xi=0.5$. The exact form of $ \bar F$
is given by
%
%
\begin{equation}\label{eqnonstd}
\bar F(x)=\frac{ 1}{ 32} W( 2x\mathrm{e}^{2} )^{ 5} x^{ -5}\qquad \mbox{for all } x\ge1,
\end{equation}
where $ W$ is the Lambert W function satisfying $ W(x)\mathrm{e}^{ W(x)}=x$ for
all $ x>0$. Observe that $ W(x)\to\infty$ as $ x\to\infty$ and $
W(x)\le\log(x)$ for $ x>1$. Furthermore,
\[
\frac{ \log(x)}{W(x) } =1+ \frac{ \log W(x)}{W(x) } \to1\qquad \mbox{as } x\to\infty,
\]
and hence $ W(x)$ is a slowly varying function. This is therefore an
example where the slowly varying term contributes significantly to $
\bar F$. That was not the case in the Pareto or the stable examples.
The result of the simulation is shown in Figure~\ref{fignonstandard}.
As expected, the choice of $ k$ plays an important role in this case,
and we see that the confidence band contains the true line when we
choose $ k=1000$ and $ \varepsilon=0.01$. Although not shown in Figure
\ref{fignonstandard}, the confidence bands perform better for smaller
values of $ k$.

\begin{figure}

\includegraphics{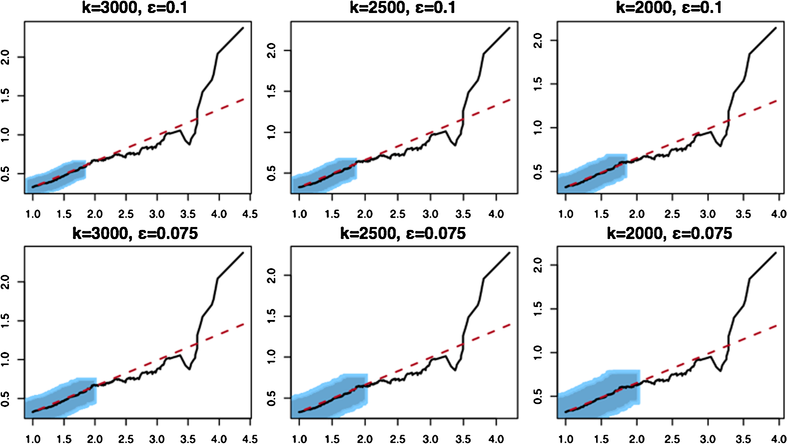}

\caption{ME plot for 50,000 i.i.d. Pareto random variables with $ \xi
=0.25$.}\label{figMEpareto}
\end{figure}

\begin{figure}

\includegraphics{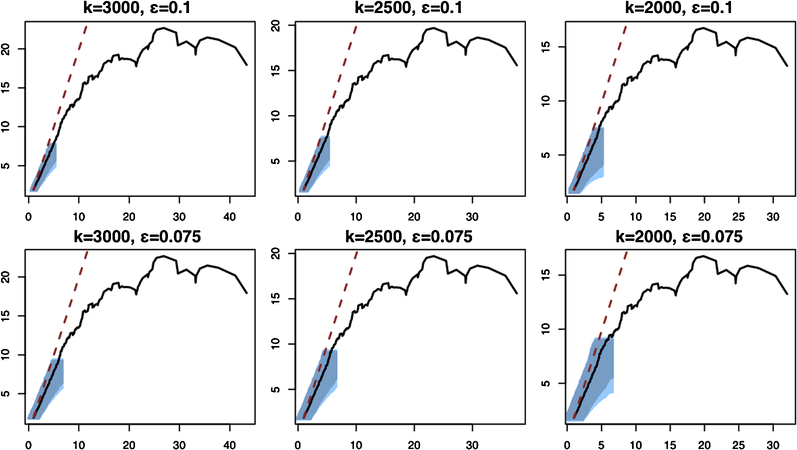}

\caption{ME plot for 50,000 i.i.d. right-skewed stable random variables
with $ \xi=2/3$.}\label{figMEstable}\vspace*{-3pt}
\end{figure}

\subsubsection{ME plots}\label{secsimulme}

Figure~\ref{figMEpareto} shows the ME plot obtained from a data
simulated from the Pareto distribution with $ \xi=0.25$. The six plots
correspond to different values of $ k$ (3000, 2500 and 2000) and $
\varepsilon$ (0.1 and 0.075). The black line is the observed ME plot,
and the brown dotted line denotes the limit in probability. Again, the
three shades of the colored bands denote the $ 99\%, 95\%$\vadjust{\goodbreak} and the $
90\%$ confidence bands for the plot, respectively. Note that the weak
limit is a functional of the Brownian bridge and depends on $ \xi$. We
estimate $ \xi$ using the Hill estimator and obtain the bounds by
simulating 10,000 paths from the weak limit.

A striking feature in all these plots is that they are close to being
linear near the bottom-left corner and become quite erratic near
top-right corner. The reason behind this phenomenon is that the
empirical ME function for high thresholds is the average of the
excesses of a small number of upper order statistics. When averaging
over few numbers, there is high variability, and therefore this part of
the plot appears very nonlinear and is uninformative. Therefore, while
obtaining confidence bands it is essential to leave out some of the
extreme order statistics. We would also like to point out that, without
the confidence bands, it would have been difficult to believe that
these plots were obtained from a distribution with tail index $ 0.25$.

A simulation of ME plot for the right skewed stable distribution with $
\xi=2/3$ is shown in Figure~\ref{figMEstable}. We use the band
described in \eqref{eqconbandMEplotxigehalf} and estimate the
quantiles using simulation. In this case we only provide the $ 95\%$
and the $ 90\%$ confidence band. The $ 99\%$ confidence band for the
stable is very large and using that is not much helpful.

\begin{figure}

\includegraphics{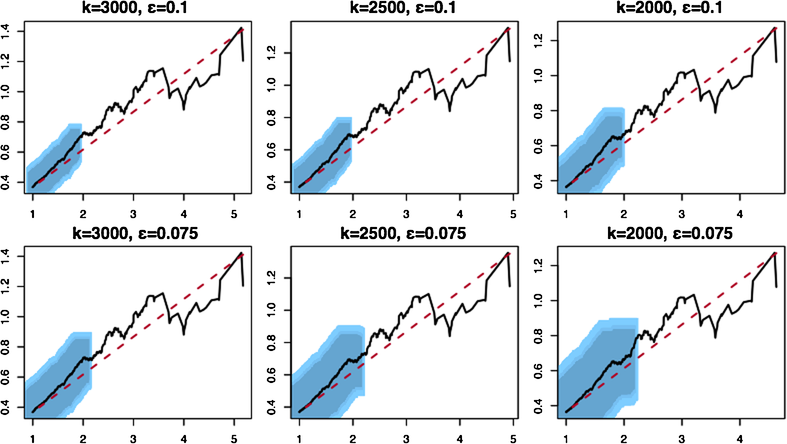}

\caption{ME plot for 50,000 i.i.d. random variables with the
distribution described in \protect\eqref{eqnonstd} ($ \xi=0.2$).}\label{figMEnonstd}
\end{figure}

The next simulation is the ME plot for a sample from the distribution
function described in~\eqref{eqnonstd}, and the result is given in
Figure~\ref{figMEnonstd}. We use the same values for $ n,k$ and~$\varepsilon$.
We see that this method of getting confidence bands works
well in these cases.

\begin{figure}

\includegraphics{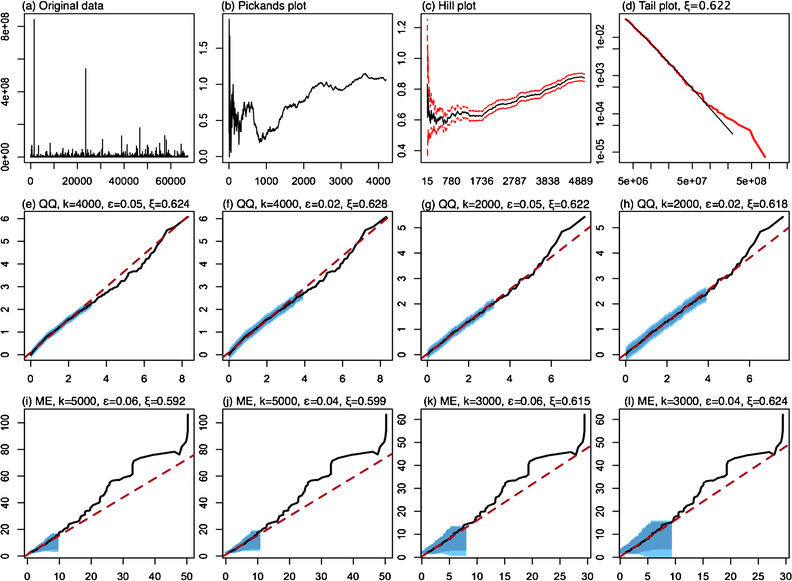}

\caption{Analysis of the internet response size data.}\label{figinternetdata}
\end{figure}

\subsection{An example with a real data}
We study a data set which contains Internet response sizes
corresponding to user requests. The sizes are thresholded to be
at least 100~KB. The data set consists of 67,287 observations and is part
of a bigger set collected in April 2000 at the University of North
Carolina at Chapel Hill.

It is often stated that file size data typically exhibits heavy tails,
and we observe that is indeed the case here. Figure \ref
{figinternetdata} shows various plots from this data set. The sample
variance is of the order of $ 10^{ 13}$ which suggests that the
variance is possibly infinite for the underlying distribution (denote
by $ F$). This would imply that if $ \bar F$ is regularly varying for
some $ \xi$, then we must have $ \xi\ge1/2$. This is suggested by both
the Pickands plot and the Hill plot (Figure~\ref{figinternetdata}(b)
and (c), resp.). The Hill plot is always above $ 1/2$ and the
Pickands is above $ 1/2$ for most of the range. But it is difficult to
get an estimate of $ \xi$ using these two tools since both plots are
highly fluctuating and hence inconclusive. We fit a GPD model with the
top 2000 order statistics using the command ``fit.GPD'' in the library
``QRMlib.'' It gives an estimate $ 0.6218$ of $ \xi$ and Figure \ref
{figinternetdata}(d) plots the estimated $ \bar F$ in the log-log
scale along with the fitted line.

We try the QQ plot with data set for $ k=4000$ and $ 2000$ (top 6\% and
3\% order statistics approximately) and with $ \varepsilon=0.05$ and $
0.02$. The plots give an estimate of around 0.62 of $ \xi$. The plots
are shown in Figure~\ref{figinternetdata}(e)--(h). The ME plots for $
k=5000, 3000$ and $ \varepsilon=0.06, 0.04$ are shown in Figure \ref
{figinternetdata}(i)--(l), and they also suggest a similar estimate
for $ \xi$.

We observe that, in this example, the different methods of
understanding the tail behavior of a data work very well, and all of
them are in agreement about the value of $ \xi$. This is not true in
many situations, and then it is hard to judge which method one should
trust. In those cases it is important to have some more knowledge about
the system from which the data was collected, and often that helps in
the understanding of the data.

\section{Conclusion}\label{secconclusion}
Plotting techniques have always been popular as diagnostic tools for
goodness-of-fit of observed data, and we believe they will remain so
because of their visual and intuitive appeal. In this paper we have
concentrated on two such tools used extensively in the extreme-value
literature. A weak law of large numbers has been shown previously for
both the QQ plots (Das and Resnick~\cite{dasresnick2008}) and ME plots
(Ghosh and Resnick~\cite{ghoshresnick2009}), considering them as random elements in an
appropriate topology. Our contribution in this paper has been to
provide distributional limits for them. In the case of QQ plots, we
have also provided an explicit expression for confidence bounds (with a
truncation to avoid the confidence bounds from blowing up) by using
these distributional results.
In the case of ME plots we have obtained distributional limits in the
cases $0<\xi<1/2$, $1/2<\xi<1$ and $\xi\ge1$ separately where the
underlying distribution $F$ is assumed to be regularly varying with
index $-1/\xi$. The case $\xi=1/2$ is still open. We have produced
confidence bounds for the ME plots in these cases by Monte Carlo
simulation, as explicit expressions for these quantities are not easy
to calculate. The explicit expressions would involve boundary-crossing
probabilities for a Brownian Bridge with nonlinear boundaries.
Boundary-crossing probabilities for Brownian motion can be approximated
using piecewise linear boundaries P\"otzelberger and Wang \cite
{potzelbergerwang2001}, but
we do not know of a nice approximation for the Brownian Bridge case;
hence we resort to simulation. We have illustrated the confidence
bounds in both the cases of QQ plots and ME plots with simulated and
real data examples in Section~\ref{secdata}. The importance of the
confidence bounds can be understood very clearly from Figure \ref
{figMEpareto}. Here we have a simulated data set of 50,000 points
from a Pareto distribution with parameter $\xi=0.25$. Just looking at
the ME plot, it is not at all obvious that this is a heavy-tailed data,
whereas when the confidence bounds with the $\varepsilon$-truncation are
drawn, the straight line with slope $\xi=0.25$ remains inside the
bounds indicating the true nature of the data.

Since we are using the limiting distribution to obtain the confidence
bounds, it is natural to ask what the rate of convergence is. We have
observed that this method works well in the simulation studies that we
have done, but we have not answered this theoretically. This is
currently a work in progress.

A standing assumption in the results we proved in this paper is that
the random variables $ X_{ n}$ are i.i.d. We believe that it is
possible to obtain similar results under a more general assumption of
stationarity and mixing; cf. Rootz\'{e}n~\cite{rootzen2009}. We intend
to look
into this further.

We should also note here that often practitioners use the
median-excess plot with the implied meaning when $\xi>1$; that is, the
mean for the distribution does not exist (Embrechts, Kl\"{u}ppelberg
and Mikosch~\cite{embrechtskluppelbergmikosch1997}), but we have not ventured into
this kind of plotting tool. We have also not looked into other kinds of
plots used in extremes, like the St\u{a}ric\u{a} plot (St\u{a}ric\u{a}
\cite{starica1999}) to determine the right $k$ number of upper order
statistics, or the Gertensgarbe and Werner plot (Gertensgarbe and
Werner~\cite{gertensgarbewerner1989}), for determining thresholds, over which a
data may be assumed to be extreme-valued, or the more popular Hill
plot, Pickands plot (Resnick~\cite{resnickbook2007}), to detect the right
value of the extreme-value parameter. Obtaining results in the same
spirit as this paper for these other varieties of plots are a part of
intended future research.

\section*{Acknowledgements}
The authors are thankful to Paul Embrechts (ETH Zurich), Sidney I.
Resnick (Cornell University) and Gennady Samorodnitsky (Cornell
University) for their detailed
comments on a draft of the paper which greatly helped in improving the
paper. The authors are also thankful for insightful comments and
suggestions from the referees and the associate editor. Bikramjit Das
was partially supported by the program IRTG/Pro*Doc. Souvik Ghosh was
partially supported by the FRAP program at Columbia University.

\def\cprime{$'$}
%
%

\printhistory

\end{document}